\documentclass{amsart}
\usepackage[latin1]{inputenc}
\usepackage[all]{xy}
\usepackage{amssymb}
\usepackage{amsmath}
\usepackage{graphicx}
\usepackage{epsfig}

\newtheorem{theorem}{Theorem}[section]
\newtheorem{rem}[theorem]{Remark}
\newtheorem{prop}[theorem]{Proposition}
\newtheorem{lemma}[theorem]{Lemma}
\newtheorem{cor}[theorem]{Corollary}
\newtheorem{defi}[theorem]{Definition}

\newcommand{\lra}{\longrightarrow}

\newcommand{\bprf}{{\it Proof.~}}

\newcommand{\ra}{\rightarrow}
\newcommand{\eprf}{\hfill $\square$ \smallskip\par}
\newcommand{\erem}{\hfill $\square$}
\newcommand{\PP}{ \mathbb{P}}

\newcommand{\R}{ \mathbb{R}}
\newcommand{\Z}{\mathbb{Z}}
\newcommand{\N}{\mathbb{N}}
\newcommand{\Q}{\mathbb{Q}}

\DeclareMathOperator{\Aut}{\rm Aut}
\DeclareMathOperator{\Sing}{\rm Sing}
\DeclareMathOperator{\rk}{\rm rank}
\DeclareMathOperator{\id}{\rm id}

\makeatletter
\def\blfootnote{\xdef\@thefnmark{}\@footnotetext}

\begin{document}

\subjclass[2010]{Primary 14J28; Secondary  14J10, 14J50.}
\keywords{Kummer surfaces, K3 surfaces, automorphisms, Enriques involutions, even sets of nodes.\\
The first author is partially supported by PRIN 2010--2011 ``Geometria delle variet\`a algebriche" and FIRB 2012 ``Moduli Spaces and Their Applications".\\
Both authors are partially supported by  GRIFGA: {\it Groupement de Recherche Italo-Fran\c cais en G\'eom\'etrie Alg\'ebrique}}

\title[Kummer surfaces, symplectic action]{Kummer surfaces and K3 surfaces with $(\Z/2\Z)^4$ symplectic action}
\author{Alice Garbagnati and Alessandra Sarti}
\address{Alice Garbagnati, Dipartimento di Matematica, Universit\`a di Milano,
  via Saldini 50, I-20133 Milano, Italia}
\email{alice.garbagnati@unimi.it}

\address{Alessandra Sarti,Universit\'e de Poitiers,
Laboratoire de Math\'ematiques et Applications,
UMR 6086 du CNRS,
T\'el\'eport 2
Boulevard Marie et Pierre Curie
BP 30179,
86962 Futuroscope Chasseneuil Cedex,
France }
\email{sarti@math.univ-poitiers.fr}
\urladdr{http://www-math.sp2mi.univ-poitiers.fr/~sarti/}

\begin{abstract}
In the first part of this paper we give a survey of classical results on Kummer surfaces with Picard number 17 from the point of view of lattice theory. We prove ampleness properties for certain divisors on Kummer surfaces and we use them to describe projective models of Kummer surfaces of $(1,d)$-polarized Abelian surfaces for $d=1,2,3$. As a consequence we prove that in these cases the N\'eron--Severi group can be generated by lines.\\ 
In the second part of the paper we use Kummer surfaces to obtain results on K3 surfaces with a symplectic action of the group $(\Z/2\Z)^4$. In particular we describe the possible N\'eron--Severi groups of the latter in the case that the Picard number  is $16$, which is the minimal possible.  We describe also the N\'eron--Severi groups of the minimal resolution of the quotient surfaces which have 15 nodes. We extend certain classical results on Kummer surfaces to these families. 
\end{abstract}

\maketitle
\section{Introduction}
Kummer surfaces are particular K3 surfaces, obtained as minimal resolutions of the quotient of an Abelian surface by an involution. They are algebraic and form a 3-dimensional family  of K3 surfaces. Kummer surfaces play a central role in the study of K3 surfaces, indeed certain results on K3 surfaces are easier to prove for Kummer surfaces (thanks to their relation with the Abelian surfaces), but can be extended to more general families of K3 surfaces: the most classical example of this is the Torelli theorem, which holds for every K3 surface. The aim of this paper is to describe some results on Kummer surfaces (some of them are classical) and to prove that these results extend to 4-dimensional families of K3 surfaces.  Every Kummer surface has the following properties: it admits the group $(\Z/2\Z)^4$ as group of automorphisms which preserves the period (these automorphisms will be called {\it symplectic}) and it is also the desingularization of the quotient of a K3 surface by the group $(\Z/2\Z)^4$ which acts preserving the period. The families of K3 surfaces with one of these properties are 4-dimensional: we study these families using the results on Kummer surfaces and we prove that several properties of the Kummer surfaces hold more in general for at least one of these families.\\

The first part of the paper  (Sections \ref{sec: generalities on Km}, \ref{sec: automorphisms on Kummer}, \ref{ampleness}, \ref{examplespic17}) is devoted to Kummer surfaces. We first recall their construction and the definition of the Shioda--Inose structure which was introduced by Morrison in \cite{morrison}. In particular we recall that every Kummer surface $Km(A)$ is the quotient of both an Abelian surface and a K3 surface by an involution (cf. Proposition \ref{shiodakummer}). Since we have these two descriptions of the same surface $Km(A)$ we obtain also two different descriptions of the N\'eron--Severi group of $Km(A)$ (see Proposition \ref{prop: NS Kummer Abelian} and Theorem \ref{prop: NS Kummer Shioda Inose}).  In Section \ref{sec: automorphisms on Kummer} we recall that every Kummer surface admits certain automorphisms, and in particular the group $(\Z/2\Z)^4$ as group of symplectic automorphisms. In Proposition \ref{prop: every Kummer Z2Z4} we show that the minimal resolution of the quotient of a Kummer surface $Km(A)$ by $(\Z/2\Z)^4$ is again $Km(A)$. This gives a third alternative description
of a Kummer surface and shows that the family of Kummer surfaces is a subfamily both of the family of K3 surfaces $X$ admitting $(\Z/2\Z)^4$ as group of symplectic automorphisms and of the family of the K3 surfaces $Y$ which are quotients of some K3 surfaces by the group $(\Z/2\Z)^4$.

The main results on Kummer surfaces are contained in Section \ref{ampleness} and applied in Section \ref{examplespic17}: Nikulin, \cite{nikulinKummer}, showed that a non empty set of disjoint smooth rational curves on a K3 surface can be the branch locus of a double cover only if it contains exactly 8  or 16 curves.  In the first case the surface which we obtain by taking the double cover and contracting the $(-1)$-curves is again a K3 surface, in the second case  the surface we obtain in the same way is an Abelian surface and the K3 surface is in fact its Kummer surface. In the sequel we call the sets of disjoint rational curves in the branch locus of a double cover {\it even set}. In \cite{projectivemodels} we studied the N\'eron--Severi group, the ampleness properties of divisors and the associated projective models of K3 surfaces which admit an even set of 8 rational curves. Here we prove similar results for K3 surfaces admitting an even set of 16 rational curves, thus for the Kummer surfaces. In Section \ref{ampleness} we prove that certain divisors on Kummer surfaces are nef, or big and nef, or ample. In Section \ref{examplespic17} we study some maps induced by the divisors considered before and we obtain projective models for the Kummer surfaces of the $(1,d)$-polarized Abelian surfaces for $d=1,2,3$. As byproduct we show that these Kummer surfaces have at least one model such that their N\'eron--Severi group is generated by lines. Several models described are already well known, but here we suggest a systematical way to produce projective models of Kummer surfaces by using lattice theory.\\

In the second part of the paper (Sections \ref{sec: K3 with Z2Z4 in general}, \ref{X}, \ref{Y}, \ref{sec: maps}, \ref{examples}) we apply the previous results on Kummer surfaces to obtain general results on K3 surfaces $X$  with symplectic action by  $(\Z/2\Z)^4$ and on the minimal resolutions $Y$ of the quotients. In Theorem \ref{reticoli}, Proposition \ref{prop: overlattices of L+MG} and Theorem \ref{preciseNS}  we describe explicitly $NS(X)$ and $NS(Y)$ and thus we describe the families of the K3 surfaces $X$ and $Y$ proving that they are 4-dimensional and specialize to the family of the Kummer surfaces.

In \cite{Keum Enriques} Keum proves that every Kummer surface admits an {\it Enriques involution} (i.e. a fixed points free involution). Here we prove that this property extends to every K3 surface $X$ admitting a symplectic action of $(\Z/2\Z)^4$ and with Picard number 16 (the minimal possible). This shows that the presence of a certain group of symplectic automorphisms on a K3 surface implies the presence of a non--symplectic involution as well. 

On the other hand certain results proved for Kummer surfaces hold also for the K3 surfaces $Y$. In Proposition \ref{prop: ampleness on Y} we prove that certain divisors on $Y$ are ample (or nef and big) as we did in Section \ref{ampleness} for Kummer surfaces. The surface $Y$ admits 15 nodes, by construction. We recall that every K3 surface with 16 nodes is in fact a Kummer surface; we prove that similarly every K3 surface which admits 15 nodes is the quotient of a K3 surface by a symplectic action of $(\Z/2\Z)^4$. This result is not trivial, indeed the following analogue is false: a K3 surface with 8 nodes is not necessarily the quotient of a K3 surface by a symplectic involution. Moreover we show in Theorem \ref{preciseNS} that K3 surfaces with 15 nodes exist for polarizations of any degree (we give some examples in Section \ref{examples}). This answers the question which is the maximal number of nodes a K3 surface with a given  polarization can have (and it does not contain further singularities). If the polarization $L$ with $L^2=2t$ has $t$ even then the maximal number is $16$ and it is attained precisely by Kummer surfaces, otherwise this maximum is $15$. 
Finally, in Section \ref{examples} we give explicit examples of the surfaces  $X$ and $Y$ and we describe their geometry.\\
We point out that some results of the paper are (partially) contained in the first author's PhD thesis, \cite{alitesi}. \\

\textbf{Acknowledgments}: We are indebted with Bert van Geemen for his support and invaluable help during the preparation of the paper. The Proposition \ref{autoku} and the Theorems \ref{preciseNS} and \ref{15nodi} are motivated by a question of Klaus Hulek and Ciro Ciliberto respectively. The study in Section \ref{examplespic17} of K3 surfaces with N\'eron--Severi group generated by lines is motivated by several discussions with Masato Kuwata.  We want to thank all of them for asking the questions and for their comments.

\section{Generalities on Kummer surfaces}\label{sec: generalities on Km}
\subsection{Kummer surfaces as quotients of Abelian surfaces}
Kummer surfaces are K3 surfaces constructed as desingularization of the quotient of an Abelian surface $A$ by an involution $\iota$. Equivalently they are K3 surfaces admitting an even set of 16 disjoint rational curves.
We recall briefly the construction: let $A$ be an Abelian surface (here we consider only the case of algebraic Kummer 
surfaces), let $\iota$ be the involution $\iota: A\lra A$, $a\mapsto -a$. Let $A/\iota$ be the quotient surface. It has sixteen singular points of type $A_1$ which are the image, under the quotient map, of the sixteen points of the set $A[2]=\{a\in A\mbox{ such that }2a=0\}$. Let $\widetilde{A/\iota}$ be the desingularization of $A/\iota$. The smooth surface $Km(A):=\widetilde{A/\iota}$ is a K3 surface. Consider the surface $\widetilde{A}$, obtained from $A$ by blowing up the points in $A[2]$. The automorphism $\iota$ on $A$ induces an automorphism $\widetilde{\iota}$ on $\widetilde{A}$ whose fixed locus are the sixteen exceptional divisors of the blow up of $A$. Hence the quotient $\widetilde{A}/\widetilde{\iota}$ is smooth. It is well known that $\widetilde{A}/\widetilde{\iota}$ is isomorphic to $Km(A)$ and that we have a commutative diagram: 
\begin{eqnarray}\label{diagkummer}
\xymatrix{
\widetilde{A}\ar[r]^-{\gamma}\ar[d]&A\ar[d]^-{\pi_A}\\
Km(A)\ar[r]&A/\iota
}
\end{eqnarray}

We observe that on $\widetilde{A}$ there are 16 exceptional curves of the blow up of the 16 points of $A[2]\subset A$. These curves are fixed by the involution $\widetilde{\iota}$ and hence are mapped to 16 rational curves on $Km(A)$. Each of these curves corresponds uniquely to a point of $A[2]$. Since $A[2]\simeq (\Z/2\Z)^4$, we denote these 16 rational curves on $Km(A)$ by $K_{a_1,a_2,a_3,a_4}$, where $(a_1,a_2,a_3,a_4)\in(\Z/2\Z)^4$. Since the points in $A[2]$ are fixed by the involution $\iota$, the exceptional curves on $\widetilde{A}$ are fixed by $\widetilde{\iota}$ and so the curves $K_{a_1,a_2,a_3,a_4}$ are the branch locus of the $2:1$ cyclic cover $\widetilde{A}\ra Km(A)$. In particular the curves $K_{a_1,a_2,a_3,a_4}$ form an even set, i.e. $\frac{1}{2}(\sum_{a_i\in\Z/2\Z}K_{a_1,a_2,a_3,a_4})\in NS(Km(A))$.
\begin{defi} (cf. \cite{nikulinKummer}). The minimal primitive sublattice of $H^2(Km(A),\Z)$ containing the 16 classes of the curves $K_{a_1,a_2,a_3,a_4}$ is called {\rm Kummer lattice} and is denoted by $K$.
\end{defi}
In \cite{nikulinKummer} it is proved that  a K3 surface $X$ is a Kummer surface if and only if the the Kummer lattice is primitively embedded in $NS(X)$.
\begin{prop}\cite[Appendix to section 5, Lemma 4]{PStorelli}
The lattice $K$ is a negative definite even lattice of rank sixteen. Its discriminant is $2^6$.
\end{prop}
\begin{rem}\label{rem: properties of K}{\rm
Here we briefly recall the properties of $K$ (these are well known and can be found e.g. in \cite{PStorelli}, \cite{bpv}, \cite{morrison}): \begin{itemize}\item[1)]
Let $W$ be a hyperplane in the affine 4-dimensional space $(\Z/2\Z)^4$, i.e. $W$ is defined by an equation of type $\sum_{i=1}^4\alpha_i a_i=\epsilon$ where $\alpha_i, \epsilon\in\{0,1\}$, and $a_i\neq 0$ for at least one $i\in\{1,2,3,4\}$. The hyperplane $W$ consists of eight points.  For every $W$, the class $\frac{1}{2}\sum_{p\in W}K_p$ is in $K$ and there are 30 classes of this kind.
\item[2)] The class $\frac{1}{2}\sum_{p\in (\Z/2\Z)^4}K_p$ is in $K$.
\item[3)] Let $W_i=\{(a_1,a_2,a_3,a_4)\in A[2]\mbox{ such that }a_i=0\}$, $i=1,2,3,4$. A set of generators (over $\Z$)  of the Kummer lattice is given by the classes: $\frac{1}{2}(\sum_{p\in(\Z/2\Z)^4}K_p)$, $\frac{1}{2}\sum_{p\in W_1}K_p$, $\frac{1}{2}\sum_{p\in W_2}K_p$, $\frac{1}{2}\sum_{p\in W_3}K_p$, $\frac{1}{2}\sum_{p\in W_4}K_p$, $K_{0,0,0,0}$, $K_{1,0,0,0}$, $K_{0,1,0,0}$, $K_{0,0,1,0}$, $K_{0,0,0,1}$, $K_{0,0,1,1}$, $K_{0,1,0,1}$, $K_{1,0,0,1}$, $K_{0,1,1,0}$, $K_{1,0,1,0}$, $K_{1,1,0,0}$.
\item[4)] The discriminant form of $K$ is isometric to the discriminant form of $U(2)^{\oplus 3}$. In particular the discriminant group is $(\Z/2\Z)^6$, there are 35 non zero elements on which the discriminant form takes value 0 and 28 non zero elements on which the discriminant form takes value 1.
\item[5)] With respect to the group of isometries of $K$ there are three orbits in the discriminant group: the orbit of zero, the orbit of the 35 non zero elements on which the discriminant form takes value 0 and the orbit of the 28 elements on which the discriminant form takes value 1.
\item[6)] Let $V$ and $V'$ be two 2-dimensional planes (they are the intersection of two hyperplanes in $(\Z/2\Z)^4$ and thus isomorphic to $(\Z/2\Z)^2$), such that $V\cap V'=\{(0,0,0,0)\} $. Denote by $V \ast V':=V\cup V'-(V\cap V')$, then the classes  $w_4:=\frac{1}{2}\sum_{p\in V}K_p$ are 35 classes in $K^{\vee}/K$ and the discriminant form on them takes value 0;  the classes $w_6:=\frac{1}{2}\sum_{p\in V\ast V'}K_p$ are in 28 classes in $K^{\vee}/K$ and  the discriminant form on them takes value 1,(see e.g. \cite[Proposition 2.1.13]{alitesi}).
\item[7)] Let $V_{i,j}=\{(0,0,0,0),\alpha_i,\alpha_j, \alpha_i+\alpha_j\}\subset (\Z/2\Z)^4$, $1\leq i,j\leq 4$ where $\alpha_1=(1,0,0,0)$, $\alpha_2=(0,1,0,0)$, $\alpha_3=(0,0,1,0)$, $\alpha_4=(0,0,0,1)$. Then $\frac{1}{2}(\sum_{p\in V_{1,2}}K_p)$, $\frac{1}{2}(\sum_{p\in V_{1,3}}K_p)$, $\frac{1}{2}(\sum_{p\in V_{1,4}}K_p)$, $\frac{1}{2}(\sum_{p\in V_{2,3}}K_p)$, $\frac{1}{2}(\sum_{p\in V_{2,4}}K_p)$, $\frac{1}{2}(\sum_{p\in V_{3,4}}K_p)$ generate the discriminant group of the Kummer lattice.
\end{itemize}}\end{rem}

Here we want to relate the N\'eron--Severi group of the Abelian surface $A$ with the N\'eron--Severi group of its Kummer surface $Km(A)$. Recall that for an abelian variety $A$ we have $H^2(A,\Z)=U^{\oplus 3}$ (see e.g. \cite[Theorem-Definition 1.5]{morrison}). 
\begin{prop}\label{prop: iota id on A}
The isometry $\iota^*$ induced by $\iota$ is the identity on  $H^2(A,\Z)$.
\end{prop}
\bprf The harmonic two forms on $A$ are $dx_i\wedge dx_j$, $i\not=j$, $i,j=1,2,3,4$ where $x_i$ are the local coordinates of $A$ viewed as the real four dimensional variety $(\R/\Z)^4$. By the definition of $\iota$ we have: $dx_i\wedge dx_j\stackrel{\iota^*}\mapsto d(-x_i)\wedge d(-x_j)=dx_i\wedge dx_j$. So $\iota$ induces the identity on $H^2(X,\R)=H^2(X,\Z)\otimes \R$ and hence on $H^2(X,\Z)$ since $H^2(X,\Z)$ is torsion free.\eprf
 Let $\widetilde{A}$ be the blow up of $A$ in the sixteen fixed points of the involution $\iota$ and let $\pi_A:A\lra A/\iota$ be the $2:1$ cover. As in \cite[Section 3]{morrison}, let $H_{\tilde A}$ be the orthogonal complement in $H^2(\tilde A,\Z)$ of the
exceptional curves and $H_{Km(A)}$ be the orthogonal complement in $H^2(Km(A),\Z)$ of the
16 $(-2)$-curves on $Km(A)$. Then  $H_{\tilde A}\cong H^2(A,\Z)$ and there are the natural maps (see \cite[Section 3]{morrison}):
$$
\pi_A^*: H_{Km(A)}\rightarrow H_{\tilde A}\cong H^2(A,\Z); \quad \pi_{A_*}:H^2(A,\Z)\cong H_{\tilde A}\rightarrow H_{Km(A)}\subset H^2(Km(A),\Z)
$$

\begin{lemma}
 We have ${\pi_A}_*(U^{\oplus 3})={\pi_A}_*(H^2(A,\Z)^{\iota^*})=H^2(A,\Z)^{\iota^*}(2)=U^{\oplus 3}(2)$ and so ${\pi_A}_*(U^{\oplus 3})=U(2)^{\oplus 3}$.
\end{lemma}
\bprf Follows from \cite[Lemma 3.1]{morrison} and Proposition \ref{prop: iota id on A}.\eprf

By this lemma we can write $\Lambda_{K3}\otimes \Q\cong H^2(Km(A),\Q)\simeq\left (U(2)^{\oplus 3}\oplus
<-2>^{\oplus 16}\right)\otimes \Q$. The lattice $U(2)^{\oplus 3}\oplus
<-2>^{\oplus 16}$ has index $2^{11}$ in  $\Lambda_{K3}\simeq U^{\oplus 3}\oplus E_8(-1)^{\oplus 2}$.

\begin{prop}\label{prop: NS Kummer Abelian}
Let $Km(A)$ be the Kummer surface associated to the Abelian surface $A$. Then the Picard number of
$Km(A)$ is $\rho(Km(A))=\rho(A)+16$, in particular $\rho(Km(A))\geq 17$.\\
The transcendental lattice of $Km(A)$ is $T_{Km(A)}=T_A(2)$. The
N\'eron--Severi group $NS(Km(A))$ is an overlattice
$\mathcal{K}_{NS(A)}'$ of $NS(A)(2)\oplus K$
and $$[NS(Km(A)):(NS(A)(2)\oplus K)]=2^{\rho(A)}.$$\end{prop}

\proof We have that $\pi_{A_*}(NS(A)\oplus T_A)=NS(A)(2)\oplus T_A(2)$ and this lattice is orthogonal to the
16 $(-2)$-classes in $H^2(Km(A),\Z)$ arising form the desingularization of $A/\iota$. \\
Since $\pi_{A_*}$ preserves the Hodge decomposition, we have $NS(A)(2)\subset NS(Km(A))$ and $T_A(2)= T_{Km(A)}$ (cf. \cite[Proposition 3.2]{morrison}).
Hence the N\'eron--Severi group of $Km(A)$ is an overlattice of finite index of $NS(A)(2)\oplus K$. In fact, $\rk NS(Km(A))=22-\rk T_A=22-(6-\rk(NS(A)))=16+\rk(NS(A))=\rk(NS(A)(2)\oplus K)$. The index of this inclusion is computed comparing the discriminant of these two lattices indeed $2^{6-\rho(A)}d(T_A)=d(T_{Km(A)})=d(NS(Km(A)))$ and
$d(NS(A)(2)\oplus K)=2^{6}2^{\rho(A)}d(NS(A))=2^{6+\rho(A)}d(T_A)$, thus $d(NS(A)(2)\oplus K)/d(NS(Km(A)))=2^{6+\rho(A)}d(T_A)/2^{6-\rho(A)}d(T_A)=2^{2\rho(A)}$ which is equal to $[NS(Km(A)):(NS(A)(2)\oplus K)]^2$ (see e.g. \cite[Ch. I, Lemma 2.1]{bpv}).\eprf

Now we will consider the generic case, i.e. the case of Kummer surfaces with Picard number 17. By Proposition \ref{prop: NS Kummer Abelian}, if $Km(A)$ has Picard number 17, then its N\'eron--Severi group is an overlattice, $\mathcal{K}_{4d}'$, of index 2 of $NS(A)(2)\oplus K\simeq \Z H\oplus K$ where $H^2=4d$, $d>0$. In the next proposition we describe the possible overlattices of $\Z H\oplus K$ with $H^{2}=4d$ and hence the possible N\'eron--Severi groups of the Kummer surfaces with Picard number 17.

\begin{theorem}\label{prop: mathcal K}
Let $Km(A)$ be a Kummer surface with Picard number 17,
let $H$ be a divisor generating  $K^{\perp}\subset NS(Km(A))$,
$H^2>0$. Let $d$ be a positive integer such that $H^2=4d$ and let
$ \mathcal{K}_{4d}:=\Z H\oplus K$. 
Then: $NS(Km(A))=\mathcal{K}_{4d}'$, where $\mathcal{K}_{4d}'$ is
generated by $\mathcal{K}_{4d}$ and by a class $(H/2,v_{4d}/2)$, with:\\
$\bullet$ $v_{4d}\in K$, $v_{4d}/2\not\in K$ and $v_{4d}/2\in K^{\vee}$ (in particular $v_{4d}\cdot K_i\in 2\mathbb{Z}$);\\ $\bullet$ $H^2\equiv -v_{4d}^2\mod8$ (in particular $v_{4d}^2\in 4\Z$).\\
The lattice $\mathcal{K}'_{4d}$ is the unique even lattice (up to
isometry) such that $[\mathcal{K}'_{4d}:\mathcal{K}_{4d}]=2$ and
$K$ is a primitive sublattice of $\mathcal{K}'_{4d}$. Hence one can assume that:\\ 
if $H^2\equiv 0\mod 8$, then $$ v_{4d}=\sum_{p\in V_{1,2}}K_p=K_{0,0,0,0}+K_{1,0,0,0}+K_{0,1,0,0}+K_{1,1,0,0};$$
if $H^2\equiv 4\mod 8$, then $$v_{4d}=\sum_{p\in (V_{1,2}\ast V_{3,4})}K_p=K_{0,0,0,1}+K_{0,0,1,0}+K_{0,0,1,1}+K_{1,0,0,0}+K_{0,1,0,0}+K_{1,1,0,0}.$$
\end{theorem}
\bprf The conditions on $v_{4d}$ to construct the lattice $\mathcal{K}_{4d}$ can be proved as in \cite[Proposition 2.1]{projectivemodels}. The uniqueness of $\mathcal{K}_{4d}'$ and the choice of $v_{4d}$ follows from the description of the orbits under the group of isometries of $K$ on the discriminant group $K^{\vee}/K$, see Remark \ref{rem: properties of K}.\eprf
\begin{rem}\label{rem: overlattice of -2^16+U(2)}{\rm (cf. \cite{alikummer}, \cite{bpv}) 1) Let $\omega_{ij}:={\pi_A}_*(\gamma^*(dx_i\wedge dx_j))$, $i<j$,
$i,j=1,2,3,4$ (we use the notation of diagram \eqref{diagkummer}). The six vectors $\omega_{i,j}$ form a basis of
$U(2)^{\oplus 3}$. The lattice generated by the Kummer lattice $K$ and by the six
classes 
$$
\begin{array}{l}
u_{ij}=\frac{1}{2}(\omega_{ij}+\sum K_{a_1,a_2,a_3,a_4})
\end{array}
$$
where the sum is over $(a_1,a_2,a_3,a_4)\in (\Z/2\Z)^4$ such that $a_i=a_j=0$, $\{i,j,h,k\}=\{1,2,3,4\}$, and $h<k$, 
is isometric to $\Lambda_{K3}$.

2) Observe that since for each $d\in\Z_{>0}$ there exist Abelian surfaces with N\'eron--Severi group isometric to $\langle 2d\rangle$, for each $d$ there exist Kummer surfaces with N\'eron--Severi group isomorphic to $\mathcal{K}_{4d}'$.}
\end{rem}

Let $\mathcal{F}_d$, $d\in\Z_{>0}$ denote the family of $\mathcal{K}_{4d}'$-polarized K3 surfaces then: 

\begin{cor}\label{cor: family Kummer} The moduli space of the Kummer surfaces has a countable number of connected components, which are the $\mathcal{F}_d$, $d\in\Z_{>0}$.\end{cor}
\proof Every Kummer surface is polarized with a lattice $\mathcal{K}_{4d}'$, for some $d$, by Propositions \ref{prop: NS Kummer Abelian} and Theorem \ref{prop: mathcal K}. On the other hand if a K3 surface is $\mathcal{K}_{4d}'$ polarized, then there exists a primitive embedding of $K$ in its N\'eron--Severi group and by \cite[Theorem 1]{nikulinKummer} it is a Kummer surface.\eprf
\begin{rem}\label{rem: divisible classes in mathcal K}
{\rm The classes of type $\frac{1}{2}(H+v_{4d}+\sum_{p\in W}K_p)$, where $H$ and $v_{4d}$ are as in Theorem \ref{prop: NS Kummer Abelian} and $W$ is a hyperplane of $(\Z/2\Z)^4$, are classes in $\mathcal{K}_{4d}'$. 
We describe this kind of classes modulo the lattice $\oplus_{p\in(\Z/2\Z)^4}\Z K_p$. We use the notation of Theorem  \ref{prop: mathcal K}.\\
If $H^2=4d \equiv 0\mod 8$, the lattice $\mathcal{K}_{4d}'$ contains:
\begin{itemize}
\item 4 classes of type $\frac{1}{2}(H-\sum_{p\in J_4}K_p)$ for certain $J_4\subset (\Z/2\Z)^4$ which contain 4 elements: these classes are $\frac{1}{2}(H+v_{4d})$ and the classes $\frac{1}{2}(H+v_{4d}+\sum_{p\in W}K_p)$ where  $W\supset \{(0,0,0,0),$ $(1,0,0,0),$ $(0,1,0,0),$  $(1,1,0,0)\}$ ;
\item 24 classes of type $\frac{1}{2}(H-\sum_{p\in J_{8}}K_p)$ for certain $J_{8}\subset (\Z/2\Z)^4$ which contain 8 elements:     
 these classes are $\frac{1}{2}(H+v_{4d}+\sum_{p\in W}K_p)$ where $W\cap \{(0,0,0,0),$ $(1,0,0,0),$ $(0,1,0,0),$  $(1,1,0,0)\}$ contains 2 elements.
\item 4 classes of type $\frac{1}{2}(H-\sum_{p\in J_{12}}K_p)$ for certain $J_{12}\subset (\Z/2\Z)^4$ which contain 12 elements: these classes are $\frac{1}{2}(H+v_{4d}+\sum_{p\in (\Z/2\Z)^4}K_p)$ and the classes $\frac{1}{2}(H+v_{4d}+\sum_{p\in W}K_p)$ where $W\cap \{(0,0,0,0),$ $(1,0,0,0),$ $(0,1,0,0),$  $(1,1,0,0)\}=\emptyset$.
\end{itemize}
If $H^2=4d \equiv 4\mod 8$, the lattice $\mathcal{K}_{4d}'$ contains:
\begin{itemize}
\item 16 classes of type $\frac{1}{2}(H-\sum_{p\in J_6}K_p)$ for certain $J_6\subset (\Z/2\Z)^4$ which contain six elements: these classes are $\frac{1}{2}(H+v_{4d})$ and the classes $\frac{1}{2}(H+v_{4d}+\sum_{p\in W}K_p)$ where $W\cap \{(1,0,0,0),(0,1,0,0),(1,1,0,0),(0,0,0,1), (0,0,1,0), (0,0,1,1)\}$ contains 4 elements;
\item 16 classes of type $\frac{1}{2}(H-\sum_{p\in J_{10}}K_p)$ for certain $J_{10}\subset (\Z/2\Z)^4$ which contain 10 elements: these classes are $\frac{1}{2}(H+v_{4d}+\sum_{p\in (\Z/2\Z)^4}K_p)$  and the classes $\frac{1}{2}(H+v_{4d}+\sum_{p\in W}K_p)$ where $W\cap$ $\{(1,0,0,0),(0,1,0,0),$ $(1,1,0,0),$ $(0,0,0,1), (0,0,1,0), (0,0,1,1)\}$ contains 2 elements.
\end{itemize}}
\end{rem}

\begin{rem}\label{rem: discriminant Km rho17}{\rm 
The discriminant group of $\mathcal{K}'_{4d}$ is generated by $(H/{4d})+\frac{1}{2}(\sum_{p\in V_{3,4}}K_p)$, $\frac{1}{2}(\sum_{p\in V_{1,3}}K_p)$, $\frac{1}{2}(\sum_{p\in V_{1,4}}K_p)$, $\frac{1}{2}(\sum_{p\in V_{2,3}}K_p)$, $\frac{1}{2}(\sum_{p\in V_{2,4}}K_p)$ if $H^2=4d\equiv 0\mod 8$ and by $(H/{4d})+\frac{1}{2}(\sum_{p\in V_{1,2}}K_p)$,
 $\frac{1}{2}(\sum_{p\in V_{1,3}}K_p)$, $\frac{1}{2}(\sum_{p\in V_{1,4}}K_p)$, $\frac{1}{2}(\sum_{p\in V_{2,3}}K_p)$, $\frac{1}{2}(\sum_{p\in V_{2,4}}K_p)$ if $H^2=4d\equiv 4\mod 8$.}\end{rem}

\subsection{Kummer surfaces as K3 surfaces with 16 nodes}\label{subsec:  K3 with 16 nodes} Let $S$ be a surface with $n$ nodes and let $\widetilde{S}$ be its minimal resolution. On $\widetilde{S}$ there are $n$ disjoint rational curves which arise from the resolution of the nodes of $S$. If $\widetilde{S}$ is a K3 surface, then $n\leq 16$,  \cite[Corollary 1]{nikulinKummer}. By \cite[Theorem 1]{nikulinKummer}, if a K3 surface admits 16 disjoint rational curves, then they form an even set and the K3 surface is in fact a Kummer surface. Conversely, as remarked in the previous section, every Kummer surface contains 16 disjoint rational curves. Thus, the Kummer surfaces are the K3 surfaces admitting the maximal numbers of disjoint rational curves or equivalently they are the K3 surfaces which admit a singular model with the maximal number of nodes.

\subsection{Kummer surfaces as quotient of K3 surfaces}\label{kummerquotient}
\begin{defi}(cf. \cite[Definition 5.1]{morrison})
 An involution $\iota$ on a $K3$ surface $Y$ is a {\rm Nikulin involution} if $\iota^*\omega=\omega$ for every $\omega\in H^{2,0}(Y)$.
\end{defi}

Every Nikulin involution has eight isolated fixed points and the minimal resolution $X$ of the quotient $Y/\iota$ is again a $K3$ surface (\cite[\S 11, Section 5]{Nikulinsymplectic}). The minimal primitive sublattice of $NS(X)$ containing the eight exceptional curves from the resolution of the singularities of $Y/\iota$ is called {\it Nikulin lattice} and it is denoted by $N$, its discriminant is $2^6$.

\begin{defi} (cf.\ \cite{vGS}) A Nikulin involution $\iota$ on a K3 surface $Y$ is a {\rm Morrison--Nikulin} involution if $\iota^*$ switches two orthogonal copies of $E_8(-1)$ embedded in $NS(Y)$.\end{defi}
By definition, if $Y$ admits a Morrison--Nikulin involution then $E_8(-1)\oplus E_8(-1)\subset NS(Y)$.
A Morrison--Nikulin involution has the following properties  (cf. \cite[Theorem 5.7 and 6.3]{morrison}): \begin{itemize}\item $T_X=T_Y(2)$; \item the lattice $N\oplus E_8(-1)$ is primitively embedded in $NS(X)$; \item the lattice $K$ is primitively embedded in $NS(X)$ and so $X$ is a Kummer surface.\end{itemize} 

\begin{defi}(cf. \cite[Definition 6.1]{morrison})
Let $Y$ be a K3 surface and $\iota$ be a Nikulin involution on $Y$. The pair $(Y,\iota)$ is a {\rm Shioda--Inose structure} if the rational quotient map 
$\pi:\xymatrix{Y\ar@{-->}[r]&X}$ is such that $X$ is a Kummer surface and $\pi_*$ induces a Hodge isometry $T_Y(2)\cong T_X$.
\end{defi}

The situation is resumed in the following diagram ($A^0$ denotes an abelian surface):
\begin{eqnarray*}
\xymatrix{
      & A^0 \ar@{->}[ld]\ar@{-->}[rd]&                   & Y   \ar@{->}[rd] \ar@{-->}[ld]&  \\
A^0/i &&X=Km(A^0) \ar[ll]\ar[rr]& &Y/\iota
}
\end{eqnarray*}
We have $T_Y\cong T_{A^0}$ by \cite[Theorem 6.3]{morrison}. 

Let $Y$ be a K3 surface and $\iota$ be a Nikulin involution on $Y$. By \cite[Theorem 5.7 and 6.3]{morrison} we conclude that
\begin{cor}
A pair $(Y,\iota)$ is a Shioda--Inose structure if and only if $\iota$ is a Morrison--Nikulin involution.
\end{cor}
For the next result see \cite[Lemma 2]{OS}.

\begin{prop}\label{shiodakummer} Every Kummer surface  is the desingularization of the quotient of a K3 surface by a Morrison--Nikulin involution, i.e.\ it is associated to a Shioda--Inose structure.\end{prop}
\begin{rem}{\rm In \cite[Lemma 5]{OS} it is proved that if $X$ is a K3 surface with $\rho(X)=20$, then each Shioda--Inose structure is induced by the same Abelian surface. This means that if $(X,\iota_1)$ and $(X,\iota_2)$ are Shioda--Inose structures
and $Y_i=Km(B_i)$ is the Kummer surface minimal resolution of $X_i/\iota_i$, $i=1,2$ then $B_1=B_2$, and so $Y_1=Y_2$. }
\end{rem}

By Proposition \ref{shiodakummer} it follows that Kummer surfaces can be defined also as  K3 surfaces which are  desingularizations of the quotients of K3 surfaces by  Morrison--Nikulin involutions. This definition leads to a different description of the N\'eron--Severi group of a Kummer surface, which we give in the following:
\begin{theorem}\label{prop: NS Kummer Shioda Inose} 
 Let $Y$ be a K3 surface admitting a Morrison--Nikulin involution $\iota$, then $\rho(Y)\geq 17$ and $NS(Y)\simeq R\oplus E_8(-1)^2$ where $R$ is an even lattice with signature $(1,\rho(Y)-17)$. Let $X$ be the desingularization of $Y/\iota$, then $NS(X)$ is an overlattice of index $2^{(\rk (R))}$ of $R(2)\oplus N\oplus E_8(-1)$.\\ 
In particular, if $\rho(Y)=17$, then: $NS(Y)\simeq \langle 2d\rangle\oplus E_8(-1)^2$, the surface $X$ is the Kummer surface of the $(1,d)$-polarized Abelian surface and the N\'eron--Severi group of $X$ is an overlattice of index 2 of $\langle 4d\rangle \oplus N\oplus E_8(-1)$.\end{theorem}
\proof 
By \cite[Theorem 6.3]{morrison} and the fact that $E_8(-1)$ is unimodular one can write
$NS(Y)=R\oplus E_8(-1)^{\oplus 2}$, with $R$ even of signature $(1,\rho(Y)-17)$. In \cite[Theorem 5.7]{morrison} it is proved that $N\oplus E_8(-1)$ is primitively embedded in $NS(X)$. Thus, arguing on the discriminant of the transcendental lattices of $Y$ and $X$ and on the lattice $R$ as in Proposition \ref{prop: NS Kummer Abelian}, one concludes the first part of the proof. For the last part of the assertion observe that the lattices $NS(Y)$ and $T_Y$ are uniquely determined by their signature and discriminant form (\cite[Theorem 2.2]{morrison}), so $T_Y=\langle -2d\rangle \oplus U^2$. By construction $T_Y(2)=T_X=T_{A^0}(2)$ so $T_{A^0}=T_Y$. This determines uniquely $NS(A^0)$, which is isometric to $\langle 2d\rangle$. Hence $A^0$ is a $(1,d)$-polarized abelian surface. 
\eprf

The overlattices $\mathcal{N}_{2d}'$ of index 2 of $\langle 2d \rangle \oplus N$ are described in \cite{projectivemodels} and, by Theorem \ref{prop: NS Kummer Shioda Inose}, we conclude that if $\rho(Y)=17$ then $NS(X)\simeq \mathcal{N}_{4d}'\oplus E_8(-1)$.

\begin{rem}{\rm Examples of Shioda--Inose structures on K3 surfaces with Picard number 17 are given e.g. in the appendix of \cite{GL}, in \cite{kumar}, \cite{vGS}, \cite{koike} and in \cite{schuett}. In all these papers the Morrison--Nikulin involutions of the Shioda--Inose structures are induced by a translation by a 2-torsion section on an elliptic fibration. In particular, in \cite{koike} all the Morrison--Nikulin involutions induced in such a way on elliptic fibrations with a finite Mordell--Weil group are classified.}
\end{rem}

\begin{rem}{\rm Proposition \ref{prop: NS Kummer Abelian} and Theorem \ref{prop: NS Kummer Shioda Inose}  give two different descriptions of the same lattice (the N\'eron--Severi group of a Kummer surface of Picard number 17): the first one is associated to the construction of the Kummer surface as quotient of an Abelian surface; the second one is associated to the construction of the same surface as quotient of another K3 surface. In general it is an open problem to pass from one description to the other, and hence to find the relation among these two constructions of a Kummer surface. However in certain cases this relation is known. In \cite{naruki} Naruki describes the N\'eron--Severi group of the Kummer surface of the Jacobian of a curve of genus 2 as in our Proposition \ref{prop: NS Kummer Abelian} and he determines a nef divisor that gives a $2:1$ map to $\mathbb{P}^2$ (we describe this map in Section \ref{jaco}). Then, 16 curves on $\mathbb{P}^2$ are constructed and it is proved that their pull backs on the Kummer surface generate the lattice $N\oplus E_8(-1)$. Similarly this relation is known if the Abelian surface is $E\times E'$, the product of two non isogenous elliptic curves $E$, $E'$. In \cite{Oguiso} the N\'eron--Severi group of $Km(E\times E')$ is described as in Proposition \ref{prop: NS Kummer Abelian}. Then the elliptic fibrations on this K3 surface are classified. In particular there exists an elliptic fibration with a fiber of type $II^*$ and two fibers of type $I_0^*$. The components of $II^*$ which do not intersect the zero section generate a lattice isometric to $E_8(-1)$ and are orthogonal to the components of $I_0^*$. The components with multiplicity 1 of the two fibers of type $I_0^*$ generate a lattice isometric to $N$ and orthogonal to the copy of $E_8(-1)$ that we have described before. Thus, one has an explicit relation between the two descriptions of the N\'eron--Severi group.}\end{rem}

\section{Automorphisms on Kummer surfaces}\label{sec: automorphisms on Kummer}
It is in general a difficult problem to describe the full automorphisms group of a given K3 surface. However for certain Kummer surfaces it is known. For example the group of automorphisms of the Kummer surface of the Jacobian of a curve of genus 2 is described in \cite{Keum} and \cite{kondo}. Similarly the group $\Aut (Km(E\times F))$ is determined in \cite{KK} in the cases: $E$ and $F$ generic and non isogenous, $E$ and $F$ generic and isogenous, $E$ and $F$ isogenous and with complex multiplication.\\ 
A different approach to the study of the automorphisms of K3 surfaces is to fix a particular group of automorphisms and to describe the families of K3 surfaces admitting such (sub)group of automorphisms. For this point of view  the following two known results (Propositions \ref{prop: every Kummer Enriques} and  \ref{prop: every Kummer Z2Z4}) which assure that every Kummer surface admits some particular automorphisms are important. Moreover, we prove also a result (Proposition \ref{autoku}) which limits the list of the admissible finite group of symplectic automorphisms on a generic Kummer surface.
\subsection{Enriques involutions on Kummer surfaces}
We recall that an {\it Enriques involution} is a fixed point free involution on a K3 surface. 

\begin{prop}\label{prop: every Kummer Enriques}{\rm(\cite[Theorem 2]{Keum Enriques})} Every Kummer surface admits an Enriques involution.\end{prop}

To prove the proposition, in \cite{Keum Enriques} the following is shown first (see \cite{Nikulinbilinear} and \cite{Ho}).

\begin{prop}\label{prop: kummer enriques embedding}{\rm(\cite[Theorem 1]{Keum Enriques})} A K3 surface admits an Enriques involution if and only if there exists a primitive embedding of the transcendental lattice of the surface in $U\oplus U(2)\oplus E_8(-2)$ such that its orthogonal complement does not contain classes with self--intersection equal to $-2$.\end{prop} 
In \cite{Keum Enriques}, the author applies the proposition to the transcendental lattice of any Kummer surface. We observe that this does not give an explicit geometric description of the Enriques involution.

\subsection{Finite groups of symplectic automorphisms on Kummer surfaces}

\begin{prop}\label{prop: every Kummer Z2Z4}{\rm(see e.g.\cite{alikummer})} Every Kummer surface $Km(A)$ admits $G=(\Z/2\Z)^4$ as group of symplectic automorphisms. These are induced by the translation by points of order two on the Abelian surface $A$ and the desingularization of $Km(A)/G$ is isomorphic to $Km(A)$ (thus every Kummer surface is also the desingularization of the quotient of a Kummer surface by $(\Z/2\Z)^4$).\end{prop}
\bprf Let $A[2]$ be the group generated by 2-torsion points. This is isomorphic with $(\Z/2\Z)^4$,
 it operates on $A$ by translation and commutes with the involution $\iota$. Hence it induces an action of $G=(\Z/2\Z)^4$ on $Km(A)$, and so on $H^2(Km(A),\Z)$. Observe that $G$ leaves the lattice $U(2)^{\oplus 3}\simeq \langle\omega_{ij}\rangle$ invariant, in fact $G$ as a group generated by translation on $A$ does not change the real two forms $dx_i\wedge dx_j$. Since $T_{Km(A)}\subset U(2)^{\oplus 3}$ the automorphisms induced on $Km(A)$ by $G$ are symplectic. Moreover since $\iota$ and $G$ commute we obtain that the surface $Km(A/A[2])$ and $\widetilde{Km(A)/G}$ are isomorphic. Finally from the exact sequence $0\ra A[2]\ra A\stackrel{\cdot 2}{\ra}A\ra 0$ we have $A/A[2]\cong A$ and so 
$\widetilde{Km(A)/G}\simeq Km(A/A[2])\simeq Km(A)$.\eprf

\begin{rem}{\rm One can also consider the quotient of $Km(A)$ by subgroups of $G=(\Z/2\Z)^4$, for example by one involution. Such an involution is induced by the translation by a point of order two. Take the Abelian surface $A\cong\R^4/\Lambda$, where $\Lambda=\langle 2e_1,e_2,e_3,e_4\rangle$ and consider the translation $t_{e_1}$ by $e_1$. Thus, $A/\langle t_{e_1} \rangle$ is the Abelian surface $B:=\R^4/\langle e_1, e_2, e_3, e_4\rangle$. 
So the desingularization of the quotient of $Km(A)$ by the automorphism induced by $t_{e_1}$ is again a Kummer surface and more precisely it is $Km(B)$. 
In particular if $NS(A)=\langle 2d\rangle$, then $NS(B)=\langle 4d\rangle$, \cite{BL}. This implies that if $NS(Km(A))\simeq \mathcal{K}'_{4d}$, then $NS(Km(B))\simeq \mathcal{K}'_{8d}$. Analogously one can consider the subgroups $G_n=(\Z/2\Z)^n\subset G$ (generated by translations), $n=1,2,3$:  if $NS(Km(A))\simeq \mathcal{K}_{4d}'$, then $NS(Km(A/G_n))\simeq\mathcal{K}'_{4\cdot 2^n\cdot d}$.}\end{rem}

\begin{prop}\label{autoku}
Let $G$ be a finite group of symplectic automorphisms of a Kummer surface $Km(A)$, where $A$ is a $(1,d)$-polarized Abelian surface and $\rho(A)=1$. Then $G$ is either a subgroup of $(\Z/2\Z)^4$, or $\Z/3\Z$ or $\Z/4\Z$.
\end{prop}
\bprf Let $G$ be a finite group acting symplectically on a K3 surface and denote by $\Omega_G$ the orthogonal complement of the $G$-invariant sublattice of the K3 lattice $\Lambda_{K3}$. An algebraic K3 surface admits the group $G$ of symplectic automorphisms if and only if $\Omega_G$ is primitively embedded in the N\'eron--Severi group of the K3 surface (cf. \cite{Nikulinsymplectic}, \cite{Ha}), hence the Picard number is greater than or equal to $\rk(\Omega_G)+1$. The list of the finite groups acting symplectically on a K3 surface and the values of $\rk(\Omega_G)$ can be found in \cite[Table 2]{Xiao} (observe that Xiao
considers the lattice generated by the exceptional curves in the minimal resolution of the quotient, he denotes its rank by $c$. This is the same as $\rk(\Omega_G)$ by \cite[Corollary 1.2]{inose}). Since we are considering Kummer surfaces such that $\rho(Km(A))=17$, if $G$ acts symplectically, then $\rk(\Omega_G)\leq 16$. This gives the following list of admissible groups $G$: $(\Z/2\Z)^i$ for $i=1,2,3,4$, $\Z/n\Z$ for $n=3,4,5,6$, $\mathcal{D}_m$ for $m=3,4,5,6$, where $\mathcal{D}_m$ is the dihedral group of order $2m$, $\Z/2\Z\times \Z/4\Z$, $(\Z/3\Z)^2$, $\Z/2\Z\times \mathcal{D}_4$, $\mathfrak{A}_{3,3}$ (see \cite{mukai} for the definition), $\mathfrak{A}_4$. We can exclude that $G$ acts symplectically on a Kummer surface for all the listed cases except $(\Z/2\Z)^i$, $i=1,2,3,4$, $\Z/3\Z$ and $\Z/4\Z$ by considering the rank and the {\it length}  
 of the lattice $\Omega_{G}$, which is the minimal number of generators of the discriminant group. For example, let us consider the case $G=\mathcal{D}_3$. The lattice $\Omega_{\mathcal{D}_3}$ is an even negative definite lattice of rank 14. Since the group $\mathcal{D}_3$ can be generated by two involutions, $\Omega_{\mathcal{D}_3}$ is the sum of two (non orthognal) copies of $\Omega_{\Z/2\Z}\simeq E_8(-2)$ and admits $\mathcal{D}_3$ as group of isometries (cf. \cite[Remark 7.9]{dihedralelliptic}). In fact $\Omega_{\mathcal{D}_3}\simeq DIH_6(14)$ where $DIH_6(14)$ is the lattice described in \cite[Section 6]{dihedral}. The discriminant group of $\Omega_{\mathcal{D}_3}\simeq DIH_6(14)$ is $(\Z/3\Z)^3\times  (\Z/6\Z)^2$, \cite[Table 8]{dihedral}. If $\mathcal{D}_3$ acts symplectically on $Km(A)$, $NS(Km(A))$ is an overlattice of finite index of $\Omega_{\mathcal{D}_3}\oplus R$ where $R$ is a lattice of rank 3. But there are no overlattices of finite index of  $\Omega_{\mathcal{D}_3}\oplus R$  with discriminant group $(\Z/2\Z)^4\times \Z/2d\Z$, which is the discriminant group of $NS(Km(A))$. Indeed, for every overlattice of finite index of $\Omega_{\mathcal{D}_3}\oplus R$ , since the rank of $R$ is 3, the discriminant group  contains at least two copies of $\Z/3\Z$.\\ 
In order to exclude all the other groups $G$ listed before, one has to know the rank and the discriminant group of $\Omega_G$: this  can be found in \cite[Proposition 5.1] {noisymplectic} if $G$ is abelian; in \cite[Propositions 7.6 and 8.1]{dihedralelliptic} if $G=\mathcal{D}_m$, $m=4,5,6$, $G=\Z/2\Z\times \mathcal{D}_4$ and $G=\mathfrak{A}_{3,3}$; in \cite[Section 4.1.1]{BG} if $G=\mathfrak{A}_4$.\eprf

\begin{rem}{\rm We can not exclude the presence of symplectic automorphisms of order 3 or 4 on a Kummer surface with Picard number 17, but we have no explicit examples of such an automorphism. It is known that there are no automorphisms of such type on $Km(A)$, if $A$ is principally polarized, cf. \cite{Keum}, \cite{kondo}.
If $Km(A)$ admits a symplectic action of $\Z/3\Z$, then $d\equiv 0\mod 3$ (this follows comparing the length of $\Omega_{\Z/3\Z}$ and of $NS(Km(A))$ as in the proof of Proposition \ref{autoku}). Moreover, the automorphism of order 3 generates an infinite group of automorphisms with any symplectic involution on $Km(A)$. Otherwise, if they generate a finite group, it has to be one of the groups listed in Proposition \ref{autoku}, but there are no groups in this list containing both an element of order 2 and one of order 3. }
\end{rem}
\subsection{Morrison-Nikulin involutions on Kummer surfaces}
Examples of certain symplectic automorphisms on a Kummer surface (the Morrison-Nikulin involutions) come from the Shioda--Inose structure. We recall that every K3 surface with Picard number at least 19 admits a Morrison--Nikulin involution. In particular this holds true for Kummer surfaces of Picard number at least 19. 
This is false for Kummer surfaces with lower Picard number. In fact since a Kummer surface with a Morrison--Nikulin involution admits also a Shioda--Inose structure as shown in Section \ref{kummerquotient} it suffices to prove:
\begin{cor}
Let $Y\cong Km(B)$ be a Kummer surface of Picard number 17 or 18, then $Y$ does not admit a Shioda--Inose structure.
\end{cor}
\bprf
If a K3 surface $Y$ admits a Shioda--Inose structure, then by Theorem \ref{prop: NS Kummer Shioda Inose} we can write 
$NS(Y)=R\oplus E_8(-1)^2$ with $R$ an even lattice of rank 1 or 2, hence the length of $NS(Y)$ satisfies $l(A_{NS(Y)})\leq 2$. 
It follows immediately that we have also $l(A_{T_Y})\leq 2$. 
Let $e_1,\ldots, e_i$, $i=5$ respectively $4$, be the generators of $T_Y$. Since $Y\cong Km(B)$ we have that $T_Y=T_B(2)$ and so the classes $e_i/2$ are independent elements of $T_Y^{\vee}/T_Y$ thus we have $2\geq l(A_{T_Y})\geq 4$, which is a contradiction. 
\eprf
In the case the Picard number is 19 we can give a more precise description of the Shioda--Inose structure: 
\begin{prop}\label{prop: ShiodaInose on Kummer} Let $Y\simeq Km(B)$ be a Kummer surface, $\rho(Y)=19$ (so $Y$ admits   a Morrison--Nikulin involution $\iota$). Let $Km(A^0)$ be the Kummer surface which is the desingularization of $Y/\iota$. Then $A^0$ is not a product of two elliptic curves.\end{prop}
\bprf If $A^0=E_1\times E_2$, $E_i$, $i=1,2$ an elliptic curve, then  the classes of $E_1$ and $E_2$ in $NS(A^0)$ span a lattice isometric to $U$. To prove that $A^0$ is not such a product  it suffices to prove that there is not a primitive embedding of $U$ in $NS(A^0)$. Assume the contrary, then $NS(A^0)=U\oplus \Z h$, so $\ell(NS(A^0))=1$. Since $Y\simeq Km(B)$ is a Kummer surface, $T_Y\simeq T_B(2)$ and thus $T_{A^0}\simeq T_B(2)$, this implies that
$1=\ell(NS(A^0))=\ell(T_{A^0})=3$ which is a contradiction.
\eprf

\section{Ampleness of divisors on Kummer surfaces}\label{ampleness}
In this section we consider projective models of Kummer surfaces with Picard number 17. The main idea is that we can check if a divisor is ample, nef, or big and nef (which is equivalent to pseudo ample) because we have a complete description of the N\'eron--Severi  group and so of the $(-2)$-curves. Hence we can apply the following criterion
 (see \cite[Proposition 3.7]{bpv}):\\

{\it Let $L$ be a divisor on a K3 surface such that $L^2\geq 0$, then it is nef if and only if $L\cdot D\geq 0$ for all  effective divisors $D$ such that $D^2=-2$.}\\

This idea was used in \cite[Proposition 3.2]{projectivemodels}, where one proves that if there exists a divisor with a negative intersection with $L$ then this divisor has self-intersection strictly less than $-2$. We refer to the description of the N\'eron--Severi group given in Proposition \ref{prop: NS Kummer Abelian}, where the N\'eron--Severi group is generated, over $\Q$, by an ample class and by 16 disjoint rational curves, which form an even set over $\Z$. Since the proofs of the next propositions are very similar to the ones given in \cite[Section 3]{projectivemodels} (where the N\'eron--Severi groups of the K3 surfaces considered are generated over $\Q$ by an ample class and by 8 disjoint rational curves forming an even set) we omit them. We denote by $\phi_L$ the map induced by the ample (or nef, or big and nef) divisor $L$ on $Km(A)$.

\begin{prop}\label{prop: pseudoample KM(A)}{\rm (cf. \cite[Proposition 3.1]{projectivemodels})}
Let $Km(A)$ be a Kummer surface such that $NS(Km(A))\simeq \mathcal{K}_{4d}'$.  Let $H$ be as in Theorem \ref{prop: mathcal K}. Then we may assume that $H$ is pseudo ample and $|H|$ has no fixed components.\end{prop}

\begin{rem}{\rm The divisor $H$ is orthogonal to all the curves of the Kummer lattice, so $\phi_H$ contracts them. The projective model associated to this divisor is an algebraic K3 surface with sixteen nodes forming an even set. More precisely $\phi_H(Km(A))$ is a model of $A/\iota$.}\end{rem}

\begin{prop}\label{ample1}{\rm (cf. \cite[Propositions 3.2 and 3.3]{projectivemodels})}
Let $Km(A)$ be  a Kummer surface such that $NS(Km(A))\simeq \mathcal{K}_{4d}'$. 
\begin{itemize}
\item If $d\geq 3$, i.e. $H^2\geq 12$, then the class $H-\frac{1}{2}(\sum_{p\in(\Z/2\Z)^4}K_p)\subset NS(Km(A))$ is an ample class. Moreover $m(H-\frac{1}{2}(\sum_{p\in(\Z/2\Z)^4}K_p))$ and $mH-\frac{1}{2}(\sum_{p\in(\Z/2\Z)^4}K_p)$ for $m\in\Z_{>0}$, are ample.
\item If $d=2$, i.e. $(H-\frac{1}{2}(\sum_{p\in(\Z/2\Z)^4}K_p))^2=0$, then $m(H-\frac{1}{2}(\sum_{p\in(\Z/2\Z)^4}K_p))$ is nef for $m\geq 1$ and $mH-\frac{1}{2}(\sum_{p\in(\Z/2\Z)^4}K_p)$ is ample for $m\geq 2$.
\end{itemize}
\end{prop}

\begin{prop}\label{ample2}{\rm (cf. \cite[Proposition 3.4]{projectivemodels})}
The divisors $H-\frac{1}{2}(\sum_{p\in(\Z/2\Z)^4}K_p)$, $mH-\frac{1}{2}(\sum_{p\in(\Z/2\Z)^4}K_p)$ and $m(H-\frac{1}{2}(\sum_{p\in(\Z/2\Z)^4}K_p))$, $m\in\Z_{>0}$, do not have fixed components for $d\geq 2$.
\end{prop}

\begin{lemma}\label{modelli}{\rm (cf. \cite[Lemma 3.1]{projectivemodels})}
The map $\phi_{H-\frac{1}{2}(\sum_{p\in(\Z/2\Z)^4}K_p)}$ is an embedding  if $H^2\geq 12$.
\end{lemma}

\begin{prop}\label{ample3} {\rm (cf. \cite[Proposition 3.5]{projectivemodels})}\text{}
\begin{itemize}
\item[1)] Let $D$ be the divisor $D=H-(K_{1}+\ldots+K_{r})$ (up to relabelling of the indices), $1\leq r \leq 16$. Then $D$ is pseudo ample for $2d>r$.
\item[2)] Let $\bar{D}=(H-K_{1}-\ldots-K_{r})/2$ with $r=4,8,12$ if $d\equiv 0\!\mod 2$ and $r=6,10$ if $d\equiv 1\!\mod 2$. Then:
\begin{itemize}
\item[$\bullet$] the divisor $\bar{D}$ is pseudo ample whenever it has positive self-intersection,
\item[$\bullet$] if $\bar{D}$ is pseudo ample then it does not have fixed components,
\item[$\bullet$] if $\bar{D}^2=0$ then the generic element in $|\bar{D}|$ is an elliptic curve.
\end{itemize}
\end{itemize}
\end{prop}

\begin{rem}{\rm In the assumptions of Lemma \ref{modelli}  the divisor $H-\frac{1}{2}(\sum_{p\in(\Z/2\Z)^4}K_p)$ defines an embedding of the surface $Km(A)$ into a projective space which sends the curves of the Kummer lattice to lines. A divisor $D$ as in Proposition \ref{ample3} defines a map from the surface $Km(A)$ to a projective space which contracts some rational curves of the even set and sends the others to conics on the image. Similarly, $\bar{D}$ defines a map from the surface $Km(A)$ to a projective space which contracts some rational curves of the even set and sends the others to lines on the image.}\end{rem}
\section{Projective models of Kummer surfaces with Picard number 17}\label{examplespic17}
Here we consider certain Kummer surfaces with Picard number 17 and we describe projective models determined
by the divisors presented in the previous section. Some of these models (but not all) are very classical.
 
\subsection{Kummer of the Jacobian of a genus 2 curve.}\label{jaco}
Let $C$ be a general curve of genus 2. It is well known that the Jacobian $J(C)$ is an Abelian surface such that $NS(J(C))=\Z L$, with $L^2=2$ and $T_{J(C)}\simeq \langle -2\rangle\oplus U\oplus U$. Hence $NS(Km(J(C)))\simeq \mathcal{K}_4'$ and $T_{Km(J(C))}\simeq \langle -4\rangle \oplus U(2)\oplus U(2)$ (see Proposition \ref{prop: NS Kummer Abelian} and Theorem \ref{prop: mathcal K}).\\
Here we want to reconsider some known projective models of $Km(J(C))$ (see \cite[Chapter 6]{GH}) using the description of the classes in the N\'eron--Severi group introduced in the previous section.

The singular quotient surface $J(C)/\iota$ is a quartic in $\mathbb{P}^3$ with sixteen nodes. For each of these nodes there exist six planes which pass through that node and each plane contains other five nodes. Each plane cuts the singular quartic surface in a conic with multiplicity 2. In this way we obtain 16 hyperplane sections which are double conics. These 16 conics are called {\it tropes}. They are the image, under the quotient map $J(C)\ra J(C)/\iota$ of different embeddings of $C$ in $J(C)$. We saw that every Kummer surface admits an Enriques involution (cf.\ Proposition \ref{prop: every Kummer Enriques}). If the Kummer surface is associated to the Jacobian of a curve a genus 2, an explicit equation of this involution on the singular model of $Km(J(C))$ in $\mathbb{P}^3$ is given in \cite[Section 3.3]{Keum Enriques}.

{\bf The polarization $H$}. The map $\phi_H$ contracts all the curves in the Kummer lattices and hence $\phi_H(Km(J(C)))$ is the singular quotient $J(C)/\iota$ in $\PP^3$. The class $H$ is the image in $NS(Km(J(C)))$ of the class generating $NS(J(C))$ (Proposition \ref{prop: NS Kummer Abelian}).
The classes corresponding to the tropes are the 16 classes (described in Remark \ref{rem: divisible classes in mathcal K}, case $4d\equiv 4\mod 8$) of the form $u_{J_6}:=\frac{1}{2}(H-\sum_{p\in J_6}K_p)$. Indeed $2u_{J_6}+\sum_{p\in J_6}K_p= H$ so they correspond to a curve in a hyperplane section with multiplicity 2; $u_{J_6}^2=-2$, so they are rational curves; $u_{J_6}\cdot H=2$, so they have degree 2. In particular the trope corresponding to the class $u_{J_6}$ passes through the nodes obtained by contracting the six curves $K_p$, where $p\in J_6$.  It is a classical result (cf. \cite[Ch. I, \S 3]{hudsonKummer}) that the rational curves of the Kummer lattice and the rational curves corresponding to the tropes in this projective models form a $16_6$ configuration of rational curves on $Km(J(C))$. This can directly checked considering the intersections between the curves $K_p$, $p\in (\Z/2\Z)^4$ and the classes $u_{J_6}$.

{\bf The polarization $H-K_{0,0,0,0}$}. Another well known model is obtained projecting the quartic surface in $\PP^3$  from a node. This gives a $2:1$ cover of $\mathbb{P}^2$ branched along six lines which are the image of the tropes passing through the node from which we are projecting. The lines are all tangent to a conic (cf. \cite[\S 1]{naruki}). 
Take the node associated to the contraction of the curve $K_{0,0,0,0}$ then the linear system associated to the projection of $J(C)/\iota$ from this node is $|H-K_{0,0,0,0}|$. The classes $u_{J_6}$  such that $(0,0,0,0)\in J_6$ are sent to lines and the curve $K_{0,0,0,0}$ is sent to a conic by the map $\phi_{H-K_{0,0,0,0}}:Km(J(C))\ra \mathbb{P}^2$. This conic is tangent to the lines which are the images of the tropes $u_{J_6}$. So the map $\phi_{H-K_{0,0,0,0}}:Km(J(C))\ra \mathbb{P}^2$ exhibits $Km(A)$ as double cover of $\mathbb{P}^2$ branched along six lines tangent to the conic $\mathcal{C}:=\phi_{H-K_{0,0,0,0}}(K_{0,0,0,0})$. The singular points of the quartic $J(C)/\iota$ which are not the center of this projection are singular points of the double cover of $\mathbb{P}^2$. So the classes $K_{a_1,a_2,a_3,a_4}$ of the Kummer lattice such that $(a_1,a_2,a_3,a_4)\neq (0,0,0,0)$ are singular points for $\phi_{H-K_{0,0,0,0}}(Km(J(C)))$ and in fact correspond to the fifteen intersection points of the six lines in the branch locus. Observe that if one fixes three of the six lines, the conic $\mathcal{C}$ is tangent to the edges of this triangle. The remaining three lines form a triangle too and the edges are tangent to the conic $\mathcal{C}$. By a classical theorem of projective plane geometry (a consequence of Steiner's theorem on generation of conics) the six vertices of the triangles are contained in another conic $\mathcal{D}$, and in fact this conic is the image of one of the tropes which do not pass through the singular point corresponding to $K_{0,0,0,0}$. This can be checked directly on $NS(Km(J(C)))$. Observe that we have in total $10$ such conics.
  
{\bf Deformation}. We observe that this model of $Km(J(C))$ exhibits the surface as a special member of the 4-dimensional family of K3 surfaces which are $2:1$ cover of $\mathbb{P}^2$ branched along six lines in general position. The covering involution induces a non-symplectic involution on $Km(J(C))$ which fixes $6$ rational curves. By Nikulin's classification of non-symplectic involutions (cf. e.g. \cite[Section 2.3]{alexeevnikulin})  the general member of the family has N\'eron--Severi group isometric to $\langle 2\rangle\oplus A_1\oplus D_{4}\oplus D_{10}$ and transcendental lattice isometric to $U(2)^{\oplus 2}\oplus \langle -2\rangle^{\oplus 2}$ which clearly contains $T_{Km(J(C))}\simeq U(2)^{\oplus 2}\oplus\langle -2\rangle$. This is a particular case of Proposition \ref{prop: Kummer subfamily X}.

{\bf The polarization $2H-\frac{1}{2}\sum_{p\in(\Z/2\Z)^4}K_p$}. We denote by $D$ this polarization. The divisor $D$ is ample by Proposition \ref{ample1}. Since $D^2=8$ the map $\phi_D$ gives a smooth projective model of $Km(J(C))$ as intersection of 3 quadrics in $\mathbb{P}^5$.  Using suitable coordinates, we can write $C$ as
$$
y^2=\prod_{i=0}^5(x-s_i)
$$
with $s_i\in\mathbb{C}$, $s_i\not=s_j$ for $i\not= j$ (it is the double cover of $\PP^1$ ramified on six points). Then by \cite[Theorem 2.5]{shioda}, 
 $\phi_D(Km(J(C)))$ has equation
\begin{eqnarray}\label{kummerjacobian1}
\left\{
\begin{array}{l}
z_0^2+z_1^2+z_2^2+z_3^2+z_4^2+z_5^2=0\\
s_0z_0^2+s_1z_1^2+s_2z_2^2+s_3z_3^2+s_4z_4^2+s_5z_5^2=0\\
s_0^2z_0^2+s_1^2z_1^2+s_2^2z_2^2+s_3^2z_3^2+s_4^2z_4^2+s_5^2z_5^2=0
\end{array}
\right.
\end{eqnarray}
in $\PP^5$.  The curves of the Kummer lattice are sent to lines by the map $\phi_D$, indeed $D\cdot K_p=1$ for each $p\in (\Z/2\Z)^4$. The image of the rational curves associated to a divisor of type $u_{J_6}$ (i.e. the curves which are tropes on the surface $\phi_H(Km(J(C)))$)  are lines: in fact one computes $D\cdot u_{J_6}=1$. So on the surface $\phi_D(Km(J(C)))$ we have 32 lines which admits a $16_6$ configuration.
Keum \cite[Lemma 3.1]{Keum} proves that the set of the tropes and the curves $K_p$, $p\in (\Z/2\Z)^4$, generate the N\'eron--Severi group (over $\Z$). Here we find the same result as a trivial application of Theorem \ref{prop: mathcal K}. Moreover we can give a geometric interpretation of this fact, indeed this implies that the N\'eron--Severi group of the surface $\phi_D(Km(J(C)))$ is generated by lines (other results about the N\'eron--Severi group of K3 surfaces generated by lines can be found e.g. in \cite{samio}). More precisely the following hold:\\
\begin{prop}\label{proplines}
The N\'eron--Severi group of the K3 surfaces which are smooth complete intersections of the three quadrics in $\mathbb{P}^5$ defined by \eqref{kummerjacobian1} is generated by lines.
\end{prop}
\bprf With the help of Theorem \ref{prop: mathcal K} we find here a set of classes generating $NS(Km(J(C)))$ which correspond to lines in the projective model of the Kummer surface $\phi_{D}(Km(J(C)))$. This set of classes  is $\mathcal{S}:=\{$ $e_1:=\frac{1}{2}(H-v_{4})$, $e_2:=\frac{1}{2}(H-K_{0,0,0,0}-K_{1,0,0,0}-K_{0,1,0,1}-K_{0,1,1,0}-K_{1,1,0,0}-K_{0,1,1,1})$, $e_3:=\frac{1}{2}(H-K_{0,0,0,0}-K_{0,1,0,0}-K_{1,1,0,0}-K_{1,0,1,0}-K_{1,0,0,1}-K_{1,0,1,1})$, $e_4:=\frac{1}{2}(H-K_{0,0,0,0}-K_{0,0,1,0}-K_{0,0,1,1}-K_{1,0,0,1}-K_{0,1,0,1}-K_{1,1,0,1})$, $e_5:=\frac{1}{2}(H-K_{0,0,0,0}-K_{0,0,0,1}-K_{0,0,1,1}-K_{1,0,1,0}-K_{1,1,1,0}-K_{0,1,1,0})$, 
$e_6:=\frac{1}{2}(H-K_{0,0,0,0}-K_{1,0,0,0}-K_{0,1,0,0}-K_{1,1,0,1}-K_{1,1,1,0}-K_{1,1,1,1})$, 
$K_{0,0,0,0}$, $K_{1,0,0,0}$, $K_{0,1,0,0}$, $K_{0,0,1,0}$, $K_{0,0,0,1}$, $K_{0,0,1,1}$, $K_{0,1,0,1}$, $K_{1,0,0,1}$, $K_{0,1,1,0}$, $K_{1,0,1,0}$, $K_{1,1,0,0}\}$. Indeed, by Theorem \ref{prop: mathcal K}, a set of generators of $NS(Km(J(C)))$ is given by $e_1$ and a set of generators of the Kummer lattice $K$ (a set of generators of $K$ is described in Remark \ref{rem: divisible classes in mathcal K}). Since for $j=2,3,4,5$ $e_j-e_1\equiv (1/2)\sum_{p\in W_{ j-1}} K_p\mod (\oplus_{p\in(\Z/2\Z)^4} \Z K_{p})$ and $e_1-e_2+e_3-e_6\equiv \frac{1}{2}\sum_{p\in(\Z/2\Z)^4}{K_p}\mod (\oplus_{p\in(\Z/2\Z)^4} \Z K_{p})$, $\mathcal{S}$ is a $\Z$-basis of $NS(Km(J(C)))$. It is immediate to check that every element of this basis has intersection 1 with $D$ and thus is sent to a line by $\phi_D$.\eprf

{\bf The nef class $H-\frac{1}{2}(H-\sum_{p\in W_i}K_p)$}. Without loss of generality we consider $i=1$ and we call this class $\bar D$. By Proposition \ref{ample3}, it defines an elliptic fibration  on $Km(J(C))$ and the eight $(-2)$-classes contained in $\bar D$ are sections of the Mordell--Weil group, the others eight 
$(-2)$-classes are components  of the reducible fibers. Observe that the class $\frac{1}{2}(H-K_{1,0,0,0}-K_{1,1,0,0}-K_{0,1,0,1}-K_{0,1,1,0}-K_{0,1,1,1}-K_{0,0,0,0})$ has self intersection $-2$, has intersection $0$ with $\bar D$ and meets
the classes $K_{1,0,0,0}$ and $K_{1,1,0,0}$ in one point. One can find easily $3$ classes more as the previous one,
so that the fibration contains $4$ fibers $I_4$. Checking in \cite[Table p. 9]{kumarell} one sees that this
is the fibration number 7 so it has no more reducible fibers and the rank of the Mordell--Weil group is $3$.

{\bf Shioda--Inose structure}.  We now describe the 3-dimensional family of K3 surfaces which admit a Shioda--Inose structure associated to  $Km(J(C))$ as described in Theorem \ref{prop: NS Kummer Shioda Inose}. It is obtained by considering K3 surfaces $X$ with $\rho(X)=17$ and with an elliptic fibration with reducible fibers $I_{10}^*+I_2$ and Mordell-Weil group equal to $\Z/2\Z$ (see Shimada's list of elliptic K3 surfaces \cite[Table 1, nr. 1343]{shimada} on the arXiv version of the paper). By using the Shioda-Tate formula (cf e.g. \cite[Corollary 1.7]{shiodamodular}) the discriminant of the N\'eron--Severi group of such a surface is $(2^2\cdot 2)/2^2$. The translation $t$ by the section of order 2 on $X$ is a Morrison--Nikulin involution, indeed it switches two orthogonal copies of $E_8(-1)\subset NS(X)$. Thus, the N\'eron--Severi group is $\langle 2d\rangle \oplus E_8(-1)\oplus E_8(-1)$, and $d=1$ because the discriminant is $2$. Hence $X$ has a Shioda--Inose structure associated to the Abelian surface $J(C)$. The desingularization of the quotient $X/t$ is the Kummer surface $Km(J(C))$ and has an elliptic fibration induced by the one on $X$, with reducible fibers $I_5^*+6I_2$ (this is the number 23 of \cite{kumarell}) and $\Z/2\Z$ as Mordell-Weil group. This Shioda--Inose structure was described in \cite[Section 5.3]{kumar}.
In Theorem \ref{prop: NS Kummer Shioda Inose} we gave a description of the N\'eron--Severi group of $Km(J(C))$ related to the Shioda--Inose structure. In particular we showed that $NS(Km(J(C))$ is an overlattice of index 2 of $\langle 4\rangle\oplus N\oplus E_8(-1)$. We denote by $Q$ the generator of $\langle 4\rangle$, by $N_i$ $i=1,\ldots, 8$ the classes of the rational curves in the Nikulin lattice $N$ and by $E_j$, $j=1,\ldots,8$ the generators of $E_8(-1)$ (we assume that $E_j$, $j=1,\ldots,7$ generate a copy of $A_7(-1)$ and $E_3\cdot E_8=1$). Then a $\Z$-basis of $NS(Km(J(C)))$ is $\{\left(Q+N_1+N_2\right)/2, N_1, \ldots, N_7, \sum_{_{i=1}^8}N_i/2, E_1,\ldots, E_8\}$. It is easy to identify a copy of $N$ and an orthogonal copy of $E_8(-1)$ in the previous elliptic fibration (the one with reducible fibers $I_5^*+I_6$); in  particular one remarks that the curves $N_i$ and $E_j$, $j=2,\ldots ,8$ are components of the reducible fibers and the curve $E_1$ can be chosen to be the zero section. This immediately gives the class of the fiber in terms of the previous basis of the N\'eron--Severi group:  $F:=Q-4E_1-7E_2-10E_3-8E_4-6E_5-4E_6-2E_7-5E_8$.

\subsection{Kummer surface of a $(1,2)$-polarized Abelian surface.}\label{subsection: Kummer 1,2}
In this section $A$ will denote always a $(1,2)$ polarized Abelian surface, and $NS(A)=\Z L$ where $L^2=4$.

{\bf The polarization $H$}.  By Proposition \ref{prop: pseudoample KM(A)} the divisor $H$ is pseudo-ample and the singular model $\phi_{H}(Km(A))$ has sixteen singular points (it is in fact $A/\iota$). Since $H^2=8$ and since by \cite[Theorem 5.2]{saintdonat} $H$ is not hyperelliptic, the K3 surface $\phi_{H}(Km(A))$  is a complete intersection of three quadrics in $\PP^5$. This model is described by Barth in \cite{barthabelian}: 
\begin{prop}{\rm \cite[Proposition 4.6]{barthabelian}} Let us consider the following quadrics:
\begin{eqnarray*}
\begin{array}{l}
Q_1=\{(\mu_1^2+\lambda_1^2)(x_1^2+x_2^2)-2\mu_1\lambda_1(x_3^2+x_4^2)+(\mu_1^2-\lambda_1^2)(x_5^2+x_6^2)=0\}\\
Q_2=\{(\mu_2^2+\lambda_2^2)(x_1^2-x_2^2)-2\mu_2\lambda_2(x_3^2-x_4^2)+(\mu_2^2-\lambda_2^2)(x_5^2-x_6^2)=0\}\\
Q_3=\{(\mu_3^2+\lambda_3^2)x_1x_2-2\mu_3\lambda_3x_3x_4+(\mu_3^2-\lambda_3^2)x_5x_6=0\}.
\end{array}
\end{eqnarray*}
Let $r=r_{1,2}r_{2,3}r_{3,1}$ where
$r_{k,j}=(\lambda_j^2\mu_k^2-\lambda_k^2\mu_j^2)(\lambda_j^2\lambda_k^2-\mu_k^2\mu_j^2)$.
If $r\neq 0$ the quadrics $Q_1$, $Q_2$, $Q_3$, generate the
ideal of an irreducible surface $Q_1\cap Q_2\cap Q_3\subset
\mathbb{P}^5$ of degree 8, which is smooth except for 16 ordinary double points and which is isomorphic to $A/\iota$.\end{prop}

The surface $A/\iota$ is then contained in each quadric of the net: $\alpha_1 Q_1+\alpha_2 Q_2+\alpha_3Q_3$, $\alpha_i\in \mathbb{C}$. We observe that the matrix $M$ associated to this net of quadrics is
a block matrix \small
\begin{eqnarray*}
\begin{array}{l}
\begin{array}{ll}
M=\left[\begin{array}{lll}B_1&0&0\\
0&B_2&0\\0&0&B_3\end{array}\right]\mbox{, where }&B_1=\left[\begin{array}{cc}\alpha_1(\mu_1^2+\lambda_1^2)+\alpha_2(\mu_2^2+\lambda_2^2)&\alpha_3(\mu_3^2+\lambda_3^2)\\
\alpha_3(\mu_3^2+\lambda_3^2)&\alpha_1(\mu_1^2+\lambda_1^2)-\alpha_2(\mu_2^2+\lambda_2^2)\end{array}\right]\end{array}\\
\begin{array}{l} B_2=\left[\begin{array}{ll}-2\alpha_1\mu_1\lambda_1-2\alpha_2\mu_2\lambda_2&-2\alpha_3\mu_3\lambda_3\\
-2\alpha_3\mu_3\lambda_3&-2\alpha_1\mu_1\lambda_1+2\alpha_2\mu_2\lambda_2\end{array}\right],\\ B_3=\left[\begin{array}{ll}\alpha_1(\mu_1^2-\lambda_1^2)+\alpha_2(\mu_2^2-\lambda_2^2)&\alpha_3(\mu_3^2-\lambda_3^2)\\
\alpha_3(\mu_3^2-\lambda_3^2)&\alpha_1(\mu_1^2-\lambda_1^2)-\alpha_2(\mu_2^2-\lambda_2^2)\end{array}\right].
\end{array}\end{array}
\end{eqnarray*}
\normalsize A singular quadric of the net is such that
$\det(M)=\det(B_1)\det(B_2)\det(B_3)=0$. One eaely check that
$\det(B_1)=\det(B_2)+\det(B_3)$. So, if $\alpha_1$, $\alpha_2$,
$\alpha_3$ are such that $\det(B_i)=\det(B_j)=0$ $i\neq j$, then
also for the third block $B_h$, $h\not= i$, $h\not= j$ one has $\det(B_h)=0$. Hence such a
choice corresponds to a quadric of rank 3. There are exactly four
possible choices of $(\alpha_1,\alpha_2,\alpha_3)\in \mathbb{C}^3$ which satisfy
the condition $\det(B_i)=0$ for $i=1,2,3$. Putting $\lambda_i=1$,
$i=1,2,3$ and
$$\begin{array}{lll}
w_1=\sqrt{(\mu_2^2-\mu_3^2)(\mu_2^2\mu_3^2-1)},\ \ &
w_2=\sqrt{(\mu_1^2-\mu_3^2)(\mu_1^2\mu_3^2-1)},\ \ &
w_3=\sqrt{(\mu_2^2-\mu_1^2)(\mu_1^2\mu_2^2-1)}
\end{array}$$
the rank 3 quadrics $S_i$ correspond to the following choices of
$(\alpha_1,\alpha_2,\alpha_3)\in \mathbb{C}^3$:
$$
\begin{array}{llllll}
S_1&\mbox{ to }&(\alpha_1,\alpha_2,\alpha_3)=(w_1,w_2,w_3)&\ \
S_2&\mbox{ to
}&(\alpha_1,\alpha_2,\alpha_3)=(w_1,w_2,-w_3)\\
S_3&\mbox{ to }&(\alpha_1,\alpha_2,\alpha_3)=(w_1,-w_2,w_3)& \ \
S_4&\mbox{ to }&(\alpha_1,\alpha_2,\alpha_3)=(w_1,-w_2,-w_3)
\end{array}
$$
Since for these choices $\det(B_i)=0$ for $i=1,2,3$, the quadrics
$S_1$, $S_2$, $S_3$, $S_4$ are of type
$(\beta_1x_1+\beta_2x_2)^2+(\beta_3x_3+\beta_4x_4)^2+(\beta_5x_5+\beta_6x_6)^2=0$, 
 the singular locus of such a quadric is the plane of $\PP^5$: 
$$\left\{\begin{array}{l}\beta_1x_1+\beta_2x_2=0\\\beta_3x_3+\beta_4x_4=0\\\beta_5x_5+\beta_6x_6=0.\end{array}\right.$$
We observe that the singular planes of $S_1$ and $S_2$ are
complementary planes in $\mathbb{P}^5$ and the same is true for
the singular planes of $S_3$ and $S_4$. Then, up to a change of
coordinates, we can assume that: $$ \begin{array}{c}S_1=y_1^2+y_2^2+y_3^2,\ \
S_2=z_1^2+z_2^2+z_3^2,\ \
S_3=(l_1y_1+m_1z_1)^2+(l_2y_2+m_2z_2)^2+(l_3y_3+m_3z_3)^2\\ A/\iota=S_1\cap S_2\cap S_3.\end{array}$$
The intersection between
$\Sing(S_1)$ and $S_2$ is a conic $C_2$. The intersection of this
conic with the hypersurface $S_3$ is made up of four points.
So $\Sing(S_1)\cap (A/\iota)=\Sing(S_1)\cap(S_1\cap S_2\cap
S_3)=\Sing(S_1)\cap S_2\cap S_3$ is made up of four points which
must be singular on $A/\iota$ (as $A/\iota$ is the complete
intersection between $S_1$, $S_2$ and $S_3$ and the points are in
$\Sing(S_1)$). These four points are four nodes of the surface
$A/\iota$. There is a complete symmetry between the four
quadrics $S_1$, $S_2$, $S_3$, $S_4$, so we have:

\begin{lemma}On each plane $\Sing(S_i)$ there are exactly four singular points of the surface $A/\iota$.\end{lemma}

Let us now consider the classes of Remark \ref{rem: divisible classes in mathcal K} described by the set 
$J_8\subset (\Z/2\Z)^4$. We call any of them $u_{J_8}$. These classes have self intersection $-2$ and they are effective. Since $u_{J_8}\cdot H=4$, they correspond to rational quartics on $A/\iota$ passing through eight nodes of the surface. Moreover, they correspond to curves with multiplicity 2, indeed $2u_{J_8}+\sum_{\in J_8} K_p$ is linearly equivalent to $H$, which is the class of the hyperplane section. The classes of these rational curves and the classes in the Kummer lattice generate the N\'eron--Severi group of $Km(A)$. These curves are in a certain sense the analogue of the tropes of $Km(J(C))$: like the tropes of $Km(J(C))$  they are rational curves obtained as special hyperplane sections of $Km(A)$ and they generate the N\'eron--Severi group of the Kummer surface together with the curves of the Kummer lattice.

{\bf The polarization $H-K_{p_1}-K_{p_2}-K_{p_3}$}. Let us choose three singular points $p_i$, $i=1,2,3$ such that $p_1$, $p_2$ are contained in $\Sing(S_1)$ and $p_3\notin \Sing(S_1)$. These three points generate a plane in $\mathbb{P}^5$. The projection of $\phi_H(Km(A))$ from this planes is associated to the linear system $H-K_{p_1}-K_{p_2}-K_{p_3}$. The map $\phi_{H-K_{p_1}-K_{p_2}-K_{p_3}}:Km(A)\ra \mathbb{P}^2$ is a $2:1$ cover of $\mathbb{P}^2$ ramified along the union of two conics and two lines. The lines are the images of two of the rational curves with classes of type $u_{J_8}$,  where $J_8$ contains $p_1,p_2,p_3\in J_8$. This description of $Km(A)$ was presented in \cite{alitesi}.

{\bf Deformation}. This model exhibits $Km(A)$ as a special member of the 6-dimensional family of K3 surfaces which are double cover of $\mathbb{P}^2$ branched along two conics and two lines. The covering involution is a non-symplectic involution fixing four rational curves. By Nikulin's classification of non-symplectic involutions (see e.g. \cite[Section 2.3]{alexeevnikulin}) it turns out that the generic member of this family of K3 surfaces has N\'eron--Severi group isometric to $\langle 2 \rangle \oplus A_1\oplus D_4^{\oplus 3}$  and transcendental lattice $U(2)^{\oplus 2}\oplus\langle -2\rangle^{\oplus 4}$ (this family is studied in details in \cite{KSTT}). The transcendental lattice  $U(2)^{\oplus 2}\oplus \langle -8\rangle $ of $Km(A)$ clearly embeds in the previous lattice. 

{\bf  The polarization $2H-\frac{1}{2}\sum_{p\in(\Z/2\Z)^4}K_p$}. We call this divisor $D$. It is ample by 
Proposition \ref{ample1}. The projective model $\phi_D(Km(A))$ is a smooth K3 surface in $\mathbb{P}^{13}$.
The curves of the Kummer lattice and the ones associated to classes of type $u_{J_8}$ are sent to lines and hence the N\'eron--Severi group of $\phi_D(Km(A))$ is generated by lines (cf. Proposition \ref{proplines}).

{\bf The nef class $\frac{1}{2}(H-\sum_{p\in J_4} K_p)$}. We call it $F$. By Proposition \ref{ample3}, it defines a map $\phi_F:Km(A)\ra \mathbb{P}^1$ which exhibits $Km(A)$ as elliptic fibration with $12$ fibers of type $I_2$ and Mordell-Weil group isomorphic to $\Z^3\oplus (\Z/2\Z)^2$. Indeed the zero section and three independent sections of infinite order are the curves $K_{a,b,c,d}$ such that $F\cdot K_{a,b,c,d}=1$. The non trivial components of the $12$ fibers of type $I_2$ are $K_{e,f,g,h}$, such that $F\cdot K_{e,f,g,h}=0$. The curves  $F+2K_{0,0,0,0}+(\sum_{p\in W_3}K_p)/2$ and $F+2K_{0,0,0,0}+(\sum_{p\in W_4}K_p)/2$ are two 2-torsion sections. This description of an elliptic fibration on $Km(A)$ follows immediately by the properties of the divisors of the N\'eron--Severi group. However a geometrical construction giving the same result is obtained considering the projection of the model of $\phi_H(Km(A))\subset\mathbb{P}^5$ from the plane $\Sing(S_1)$. The image of this projection lies in the complementary plane $\Sing(S_2)$ and is a conic $C$. Let $p$ be a point of $C$ and let $\mathbb{P}^3_p$ be the space generated by $\Sing(S_1)$ and by $p$. The fiber over $p$ is $S_2\cap S_3\cap \mathbb{P}^3_p$. The fiber over a generic point of $C$ is an elliptic curve (the intersection of two quadric in $\mathbb{P}^3$). There are 12 points in $C$, corresponding to the 12 singular points of $\phi_H(Km(A))$ which are not on the plane $\Sing(S_1)$, such that the fibers over these points are singular and in fact of type $I_2$. A geometrical description of this elliptic fibration is provided also in \cite{mehran}, where it is obtained as double cover of an elliptic fibration on $Km(J(C))$.

{\bf Shioda--Inose structure}.  We now describe the 3-dimensional family of K3 surfaces which admit a Shioda--Inose structure associated to  $Km(A)$ as described in Theorem \ref{prop: NS Kummer Shioda Inose}. It is obtained using results of \cite[Section 4.6]{vGS}: consider the K3 surface $X$ with $\rho(X)=17$ and admitting an elliptic fibration with fibers $I_{16}+8I_1$ and Mordell-Weil group isometric to $\Z/2\Z$. By \cite[Proposition 4.7]{vGS} the discriminant of $NS(X)$ is $4$ and the translation $t$ by the 2-torsion section is a Morrison--Nikulin involution. Thus, the desingularization of $X/t$ is a Kummer surface, which is in fact $Km(A)$ by Theorem \ref{prop: NS Kummer Shioda Inose}. The elliptic fibration induced on $Km(A)$ has $I_8+8I_2$ singular fibers and Mordell--Weil group $(\Z/2\Z)^2$. Using the curves contained in the elliptic fibration one can easily identify the sublattice $N\oplus E_8(-1)$ of  $NS(Km(A))$: the lattice $N$ contains the 8 non trivial components of the 8 fibers of type $I_2$ and the lattice $E_8(-1)$ is generated by 7 components of the fiber of type $I_8$ and by the zero section.\\
As in the case of the Jacobian of a curve of genus 2, we give a $\Z$-basis of the N\'eron--Severi group of $Km(A)$ related to the Shioda--Inose structure and we identify the class of the fiber of this fibration: with the previous notation a $\Z$-basis is given by $\{\langle(Q+N_1+N_2+N_3+N_4\rangle)/2, N_1, \ldots, N_7, \sum_{_{i=1}^8}N_i/2, E_1,\ldots, E_8\}$, where $Q^2=8$ and $Q$ is orthogonal to $N\oplus E_8(-1)$; the class of the fiber in terms of the previous basis of the N\'eron--Severi group is  $F:=Q-5E_1-10E_2-15E_3-12E_4-9E_5-6E_6-3E_7-8E_8$.  
\subsection{Kummer surface of a $(1,3)$ polarized Abelian surface.}
Let $A$ be a $(1,3)$ polarized Abelian surface, then $NS(A)=\Z L$, $L^2=6$.

{\bf The polarization $H$}. The model of the singular quotient $A/\iota$ is associated to the divisor $H$ in $NS(Km(A))$ with $H^2=12$. By Proposition \ref{prop: pseudoample KM(A)} and  \cite[Theorem 5.2]{saintdonat} this model is a singular K3 surface in $\mathbb{P}^7$. Let us now consider the 16 classes of Remark \ref{rem: divisible classes in mathcal K} associated to the set $J_{10}\subset (\Z/2\Z)^4$. We call any of them $u_{J_{10}}$. 
They are $(-2)$-classes (see Remark \ref{rem: divisible classes in mathcal K}) and are sent to rational curves of degree $6$ on $\phi_H(Km(A))$.

{\bf The polarization $H-\frac{1}{2}(\sum_{p\in(\Z/2\Z)^4}K_p)$}. We call it $D$. It is ample by Proposition \ref{ample1} and since $D^2=4$, the surface $\phi_D(Km(A))$ is a smooth quartic in $\mathbb{P}^3.$ The curves of the Kummer lattice and the curves associated to $u_{J_{10}}$ are sent to lines. Since the classes of the curves in the Kummer lattice and the classes $u_{J_{10}}$ generate the N\'eron--Severi group of $Km(A)$, the N\'eron--Severi group of $\phi_D(Km(A))$ is generated by lines (cf. Proposition \ref{proplines}). 

{\bf The polarization $H-K_{0,0,1,0}-K_{0,0,1,1}-K_{1,0,0,0}- K_{0,1,0,0}-K_{0,0,1,1}$}. It defines a $2:1$ map from $Km(A)$ to $\mathbb{P}^2$, since 11 curves $K_p$ are contracted the branch locus is a reducible
sextic with 11 nodes.

{\bf Deformation}. The generic K3 surface double cover of $\PP^2$ branched on a reducible
sextic with 11 nodes lies in a 8-dimensional family and has transcendental lattice equal to $U(2)^{\oplus 2}\oplus \langle -2\rangle^{\oplus 6}$, see \cite[Section 2.3]{alexeevnikulin}. Clearly the transcendental lattice $U(2)^{\oplus 2}\oplus \langle -12\rangle$ 
can be primitively embedded in $U(2)^{\oplus 2}\oplus \langle -2\rangle^{\oplus 6}$, so the family of Kummer surfaces of a $(1,3)$-polarized Abelian surface is a special 3-dimensional subfamily.

{\bf The nef class $\frac{1}{2}(H-\sum_{p\in J_6} K_p)$}. We call it $F$. By Proposition \ref{ample3} it defines an elliptic fibration $Km(A)\ra \mathbb{P}^1$ with 10 fibers of type $I_2$: the components of these fibers not meeting the zero section are the curves $K_{a,b,c,d}$ of the Kummer lattice such that $F\cdot K_{a,b,c,d}=0$. The Mordell--Weil group is $\Z^5$ and the curves $K_{e,f,g,h}$ such that $F\cdot K_{e,f,g,h}=1$ are the zero section and 5 sections of infinite order (but they are not the $\Z$-generators of the Mordell--Weil group).

{\bf Shioda--Inose structure}.We now describe the 3-dimensional family of K3 surfaces which admit a Shioda--Inose structure associated to  $Km(A)$ as described in Theorem \ref{prop: NS Kummer Shioda Inose}. It was already described independently in \cite[Remark 3.3.1 (Section 3.3)]{alitesi} and \cite[Section 3.1]{koike}. Let us consider the K3 surfaces $X$ with $\rho(X)=17$ and with an elliptic fibration with reducible fibers $I_6^*+I_6$ and Mordell--Weil group $\Z/2\Z$ (as in the arXiv version of the paper \cite[Table 1, nr. 1357]{shimada}). The translation $t$ by the 2-torsion section is a Morrison--Nikulin involution (in fact it is immediate to check that it switches two orthogonal copies of $E_8(-1)\subset NS(X)$)  and hence the desingularization of the quotient $X/t$ is a Kummer surface. The latter admits an elliptic fibration induced by the one on $X$, with reducible fibers $I_3^*+I_3+6I_2$ and a 2-torsion section. By the Shioda-Tate formula (see e.g. \cite[Corollary 1.7]{shiodamodular}) the discriminant of the N\'eron--Severi group of such an elliptic fibration is $(4\cdot 3\cdot 2^6)/2^2$ and thus this Kummer surface is the Kummer surface of a $(1,3)$-polarized Abelian surface. As in the case of the Jacobian of a curve of genus 2, we give a $\Z$-basis of the N\'eron--Severi group of $Km(A)$ related to the Shioda--Inose structure and we can identify the class of the fiber of this fibration: the 8 curves $N_i$ are the 6 non trivial components of each fiber of type $I_2$ and 2 non trivial components of $I_3^*$ with multiplicity 1; the curves $E_i$ are the zero section, two components of $I_3$ and five components of $I_3^*$.
With the previous notation a $\Z$-basis is given by $\{\langle(Q+N_1+N_2\rangle)/2, N_1, \ldots N_7, \sum_{_{i=1}^8}N_i/2, E_1,\ldots, E_8\}$, where $Q^2=12$ and $Q$ is orthogonal to $N\oplus E_8(-1)$; the class of the fiber in terms of this basis of the N\'eron--Severi group is $F:=Q-6E_1-12E_2-18E_3-15E_4-12E_5-8E_6-4E_7-9E_8$.

\section{K3 surfaces with symplectic action of the group $(\Z/2\Z)^4$ and their quotients}\label{sec: K3 with Z2Z4 in general}
In the following sections we study two 4-dimensional families of K3 surfaces  that contain subfamilies of Kummer surfaces. Indeed, we have seen  that every Kummer surface admits a symplectic action of the group $(\Z/2\Z)^4$ (Proposition \ref{prop: every Kummer Z2Z4}), but the moduli space of K3 surfaces with symplectic action by $(\Z/2\Z)^4$ has dimension 4 and thus the Kummer surfaces are a 3-dimensional subfamily. We will also study the family of K3 surfaces obtained as desingularization of the quotient of a K3 surface by the group $(\Z/2\Z)^4$ acting symplectically on it. By Proposition \ref{prop: every Kummer Z2Z4} this family also contains the 3-dimensional family of Kummer surfaces.\\

Let $G=(\Z/2\Z)^{4}$ be a group of symplectic automorphisms on a K3 surface $X$. We observe that $G$ contains $(2^4-1)=15$
symplectic involutions so we have $8\cdot 15=120$ distinct points with
non trivial stabilizer group on X, and these are all the points with a non trivial stabilizer on $X$ (cf.
\cite[Section 5]{Nikulinsymplectic}). Moreover we have a commutative diagram:

\begin{eqnarray}\label{diagrammone}
\begin{array}{ccc} \widetilde{X}&\stackrel{\beta}{\lra}&X\\
\pi \downarrow&&\downarrow \pi'\\
Y&\stackrel{\tilde{\beta}}\lra& \bar{Y},
\end{array}
\end{eqnarray}
where $\bar{Y}$ is the quotient of $X$ by $G$, $\widetilde{X}$ is the blow up of $X$ at the 120 points with non trivial stabilizer (hence it contains $120$ $(-1)$-curves) and $Y$ is the minimal resolution of the quotient $\bar{Y}$ and simultaneously the quotient of  $\widetilde{X}$ by the induced action. Observe that $Y$ contains $15$ $(-2)$-curves coming from the resolution of the singularities. In fact each fixed point on $X$ has a $G$-orbit of length 8. In particular the rank of the N\'eron--Severi group of $Y$ is at least $15$ and in fact $16$ if $X$, and so $Y$, is algebraic. In particular, since  by \cite[Corollary 1.2]{inose} $\rk NS(X)$=$\rk NS(Y)$, a K3 surface with a symplectic action of $(\Z/2\Z)^{4}$ has at least Picard number $15$ ($16$ if it is algebraic).
Finally $\pi$ is $16:1$ outside the branch locus.

\section{K3 surfaces with symplectic action of $(\Z/2\Z)^4$}\label{X}
In this section we analyze the K3 surface $X$ admitting a symplectic action of $(\Z/2\Z)^4$, in particular we identify the possible N\'eron--Severi groups of such a K3 surface if the Picard number is 16, which is the minimum possible for an algebraic K3 surface with this property. This allows us to describe the families of such K3 surfaces (cf. Corollary \ref{corollary: family X}) and to prove that every K3 surfaces admitting $(\Z/2\Z)^4$ as group of symplectic automorphisms also admits an Enriques involution: this generalizes the similar result  for Kummer surfaces given in  Proposition \ref{prop: every Kummer Enriques}.

\subsection{The N\'eron--Severi group of $X$}
\begin{theorem}{\rm (cf.\ \cite{alitesi})}\label{reticoli}
Let $X$ be an algebraic K3 surface with a symplectic action of $(\Z/2\Z)^4$ and let $\Omega_{(\Z/2\Z)^4}^{\perp}=<-8>\oplus U(2)^{\oplus 3}$ be the invariant
lattice $H^2(X,\Z)^{(\Z/2\Z)^4}$. Then
$\rho(X)\geq 16$. If $\rho(X)=16$, denote by $L$ a generator of $(\Omega_{(\Z/2\Z)^4})^{\perp}\cap NS(X)$
with $L^2=2d>0$. Let
$$
\mathcal{L}^{2d}_{(\Z/2\Z)^4}:=\Z L\oplus \Omega_{(\Z/2\Z)^4}\subset NS(X).
$$
Denote by $\mathcal{L}'^{2d}_{(\Z/2\Z)^4,r}$ an overlattice of $\mathcal{L}^{2d}_{(\Z/2\Z)^4}$ of index $r$. Then there are the following possibilities for $d$, $r$ and $L$.
\begin{itemize}
\item[1)] If $d\equiv 0 \mod 2$ and $d\not \equiv 4 \mod 8$,  then $r=2$,  $L=w_1:=(0,1,t,0,0,0,0)\in\Omega_{(\Z/2\Z)^4}^{\perp}$ and  $L^2=w_1^2=4t$.\\
\item[2)] If $d\equiv 4 \mod 8$ and $d\not\equiv -4\mod 32$, then:\\ either $r=2$, $L=w_1:=(0,1,t,0,0,0,0)\in\Omega_{(\Z/2\Z)^4}^{\perp}$ and $L^2=w_1^2=4t$,\\ or $ r=4 $, $L=w_2:=(1,2,2s,0,0,0,0)\in\Omega_{(\Z/2\Z)^4}^{\perp}$ and $L^2=w_2^2=8(2s-1)$.\\
\item[3)] If $d\equiv -4\mod 32$ then:\\ either $r=2$, $L=w_1:=(0,1,t,0,0,0,0)\in\Omega_{(\Z/2\Z)^4}^{\perp}$ and $L^2=w_1^2=4t$,\\ or $ r=4 $, $L=w_2:=(1,2,2s,0,0,0,0)\in\Omega_{(\Z/2\Z)^4}^{\perp}$ and $L^2=w_2^2=8(2s-1)$,\\ or $ r=8 $, $L=w_3:=(1,4,4u,0,0,0,0)\in\Omega_{(\Z/2\Z)^4}^{\perp}$ and $L^2=w_3^2=8(8u-1)$.
\end{itemize}
If $NS(X)$ is an overlattice of $\Z w_1\oplus \Omega_{(\Z/2\Z)^4}$, then $T_X\simeq \langle-8\rangle\oplus\langle-4t\rangle\oplus U(2)^{\oplus 2}$;\\
If $NS(X)$ is an overlattice of $\Z w_2\oplus \Omega_{(\Z/2\Z)^4}$, then $T_X\simeq \left[\begin{array}{rr}-8&4\\4&-4s\end{array}\right]\oplus U(2)^{\oplus 2}$;\\
If $NS(X)$ is an overlattice of $\Z w_3\oplus \Omega_{(\Z/2\Z)^4}$, then $T_X\simeq \left[\begin{array}{rr}-8&2\\2&-4u\end{array}\right]\oplus U(2)^{\oplus 2}$.
\end{theorem}
\bprf
Since $\Omega_{(\Z/2\Z)^4}\subset NS(X)$  and $X$ is algebraic we have $\rho(X)\geq 16$.
The proof of the unicity of the possible overlattices of
$\mathcal{L}^{2d}_{(\Z/2\Z)^4}$ is based on the following idea.
Let us consider
the lattice orthogonal to $\Omega_{(\Z/2\Z)^4}$ in $\Lambda_{K3}$. For
each element $s(=L)\in\Omega_{(\Z/2\Z)^4}^{\perp}$ in a different orbit under isometries of
$\Omega_{(\Z/2\Z)^4}^{\perp}$, we can consider the
lattice $\Z s\oplus \Omega_{(\Z/2\Z)^4}$. To compute the index of the the overlattice 
$R(=NS(X))$ of $\Z s\oplus \Omega_{(\Z/2\Z)^4}$  which is primitively embedded in $\Lambda_{K3}$, we consider 
the lattice $R^{\perp}=s^{\perp}\cap \Omega_{(\Z/2\Z)^4}^{\perp}=(\Z s\oplus\Omega_{(\Z/2\Z)^4)})^{\perp}\subset {\Lambda_{K3}}$
(which is isometric to $T_X$). We compute then the discriminant group of $R^{\perp}$ to get the discriminant group
of $R$ and so we get the index $r$ of $\Z s\oplus \Omega_{(\Z/2\Z)^4}$ in $R(=NS(X))$. Recall that 
$$\Omega_{(\Z/2\Z)^4}^{\perp}\simeq \langle
-8\rangle\oplus U(2)^3\simeq (\langle -4\rangle\oplus U^3)(2).
$$
The orbits of elements by isometries of this lattice are determined by the orbits of elements by isometries 
of the lattice $\langle -4\rangle\oplus U^3$. In the
next sections we investigate them, then the proof of the theorem follows from the results
of Section \ref{section: the lattice 2d+u+u}. 
We remark moreover that under our assumptions two overlattices $R_i\supset \Z w_i\oplus \Omega_{(\Z/2\Z)^4}$
and $R_j\supset \Z w_j\oplus \Omega_{(\Z/2\Z)^4}$, $i\neq j$, cannot be isometric in $\Lambda_{K3}$ since their orthogonal
complements $R_i^\perp$ and $R_j^{\perp}$ are different. These are determined in Proposition \ref{prop: possible orbit}  below and they are the transcendental lattices $T_X$ in our statement. 
\eprf
\subsection{The lattice $\langle -2d\rangle\oplus U\oplus
U$}\label{section: the lattice 2d+u+u}
\begin{lemma}\label{lemma: u+u orbit in u}\index{Kummer surface! rank 18}
Let $(a_1,a_2,a_3,a_4)$ be a vector in the lattice $U\oplus U$. There exists an
isometry which sends the vector $(a_1,a_2,a_3,a_4)$ to the vector
$(d,de,0,0)$. In particular the vector $(a_1,a_2,0,0)$ can be sent
to $(d,de,0,0)$ where $d=\gcd(a_1,a_2)$ and $d^2e=a_1a_2$.
\end{lemma}
\bprf The lattice $U\oplus U$ is isometric to the lattice
$\{M(2,\Z), 2\det\}$ of the square matrices of dimension two
with bilinear form induced by the quadratic form given by the determinant multiplied by $2$. Explicitly the
isometry can be written as
$$
U\oplus U\lra M(2,\Z), \,\,\,
\left(\left(\begin{array}{l}a_1\\a_2\end{array}\right),\left(\begin{array}{l}a_3\\a_4\end{array}\right)\right)\mapsto\left[\begin{array}{rr}a_1&-a_3\\
a_4&a_2\end{array}\right].$$ It is well known  that under the
action of the orthogonal group $O(M(2,\Z))$ each matrix of
$M(2,\Z)$ can be sent in a diagonal matrix with diagonal
$(d_1,d_2)$, $d_1|d_2$ (this is the Smith Normal Form). Thus the lemma follows.\erem

\begin{lemma}\label{lemma: isometry on 4+U+U}
There exists an isometry which sends the
primitive vector $(a_0,a_1,a_2,a_3,a_4)\in T_{2d}:=\langle -2d\rangle\oplus
U\oplus U$, to a primitive vector
$(a,d,de,0,0)\in \langle -2d\rangle\oplus
U\oplus U$.
\end{lemma}
\bprf The primitive vector $(a_0,a_1,a_2,a_3,a_4)$ is sent to a primitive vector by any isometry.  By Lemma \ref{lemma: u+u orbit in u} there exists an isometry sending $(a_1,a_2,a_3,a_4)\in U\oplus U$ to $(d,de,0,0)\in U\oplus U$, thus there exists an isometry sending $(a_0,a_1,a_2,a_3,a_4)$ to $(a_0,d,de,0,0)$ and $(a_0,d,de,0,0)$ is primitive.  \eprf

The previous lemma allows us to restrict our attention to the vectors in the lattice
$\mathcal{A}_{2d}:=\langle -2d\rangle\oplus U$.

\begin{lemma}\label{lemma 0,1,c} There exists an isometry of $\mathcal{A}_{2d}$ which sends the vector
$(a,1,c)$, to the vector $(0,1,r)$, where $2c-2da^2=2r$.
\end{lemma}
\bprf 
First we observe that $(a,1,c)\cdot (a,1,c)=(0,1,r)\cdot
(0,1,r)=2r$. Let $R_v$ denote the reflection with respect to $v=(1,0,d)$, then for 
$w=(x,y,z)$ we have
\begin{eqnarray*}
R_v(w)=w-2\frac{w\cdot v}{v\cdot
v}v=\left(\begin{array}{c}-x+y\\y\\-2dx+dy+z\end{array}\right).\end{eqnarray*}
If $a>0$ we apply the reflection $R_v$ to $(a,1,c)$, ($v=(1,0,d)$):
$$R_v\left(\begin{array}{c} a\\1\\c\end{array}\right)=\left(\begin{array}{c} 1-a\\1\\-2da+d+c\end{array}\right).$$
Let $D$ be the isometry of $A_{2d}$,
$$D=\left[\begin{array}{rrr}-1&0&0\\0&1&0\\0&0&1\end{array}\right].$$
Then $$D\circ
R_v\left(\begin{array}{c}a\\1\\c\end{array}\right)=\left(\begin{array}{c}
a-1\\1\\-2da+d+c\end{array}\right).$$ Applying $a$ times the
isometry $D\circ R_v$ we obtain
$$(D\circ
R_v)^{a}\left(\begin{array}{lll}a\\1\\c\end{array}\right)=\left(\begin{array}{lll}0\\1\\2r\end{array}\right).$$
\erem

\begin{lemma}\label{lemma 1,2,c} There exists an isometry of
$\mathcal{A}_{2d}$ which sends a vector $q_2:=(wh\pm j,w,wt)$, with
$t,h\in\Z$, $w,j\in\N$, $0<j\leq\llcorner d/2\lrcorner$ to the
vector $p_2:=(j,w,s)$, where $s=-dwh^2\mp 2dhj+wt$.
\end{lemma}
\bprf Without lost of generality we can assume $h>0$ (if $h\leq 0$, it is
 sufficient to consider the action of $D$). Let us apply the isometry $D\circ R_v$ to the vector $q_2$:
$$(D\circ R_v)\left(\begin{array}{c}wh\pm j\\w\\wt\end{array}\right)=\left(\begin{array}{c} w(h-1)\pm j\\w\\-2d(wh\pm j)+dw+wt\end{array}\right).$$
As in the previous proof, applying $D\circ R_v$ decreases the
first component and the second remains the same. Applying
$h$-times the isometry to $q_2$, we obtain that the first
component is $j$ or $-j$. In the second case we apply again the
isometry $D$, and so in both situations we obtain $p_2$.\erem

\begin{lemma}\label{lemma: (n,b,bf,0,0)} Let $p$ be a prime
number. Let us consider the lattice $T_{2p}=\langle
-2p\rangle\oplus U\oplus U$. There exists an isometry of $T_{2p}$
which sends the vector $q:=(n,b,bf,0,0)$, $b\in\Z_{>0}$, 
$n\in\N$, $\gcd(n,b)=1$ in one of the following vectors:\\
$\bullet$ $v_{1}=(0,1,r,0,0)$ where $2b^2f-2pn^2=2r$;\\ $\bullet$
$v_2=(1,2,2s,0,0)$, where $2b^2f-2pn^2=8s-2p$;\\ $\bullet$
$v_p=(l,p,pt,0,0)$, where $2b^2f-2pn^2=2pt-2pl^2$, $0<l\leq
\llcorner p/2\lrcorner$;\\
$\bullet$ $v_{2p}=(j,2p,2pu,0,0)$, where $2b^2f-2pn^2=8p^2u-2pj^2$,
$0<j< p$, $j\equiv 1\mod 2$.
\end{lemma} \bprf We can assume $n\in\N$ and $b>0$ (if it is not
the case it suffices to consider the action of $-\id$ and of $D$). Let us
consider the reflection $R_v$, associated to the vector
$v=(1,0,p,0,0)$. We have
$$R_v\left(\begin{array}{c}n\\b\\bf\\0\\0\end{array}\right)=\left(\begin{array}{c}-n+b\\b\\-2pn+pb+bf\\0\\0\end{array}\right).$$
Again we can change the sign of the first component and we obtain
$(b-n,b,-2pn+pb+bf,0,0)$. By Lemma \ref{lemma: u+u orbit in u}
this vector can be transformed in $(b-n,b_1,b_1f_1,0,0)$, where
$\gcd(b,-2pn+pb+bf)=b_1$. Then $b_1\leq b:=b_0$. We apply now Lemma \ref{lemma: u+u orbit in u}
to the vector $(n_1, b_1, b_1f_1,0,0)$ with $n_1:=|b-n|>0$ (eventually change the sign of $b_n$
by using the matrix $D$). 
The second component of the vector
$b_1$ is a positive number, so after a finite number of
transformations there exists $\eta$ such that $b_\eta=b_{\eta+1}$ and $\gcd (n_\eta,b_\eta)=1$. 
Since $b_\eta|(pb_\eta+b_\eta f)$ and $\gcd (n_\eta,b_\eta)=1$ (recall that the image of a primitive vector
by an isometry is again primitive) $b_\eta=b_{\eta+1}$ if and only if $b_{\eta}$ divides $2p$, i.e. if
$b_\eta=1,2,p,2p$. Moreover $\gcd(b_\eta-n_\eta,b_\eta)=1$. So by Lemma
\ref{lemma: isometry on 4+U+U} and applying eventually the transformation $D$ to get the first
component of the vector positive,  after a finite number of
transformations we obtain that $q$ is isometric to one of the
vectors $(a,1,f',0,0)$, $(2k+1,2,2f',0,0)$, $(ph\pm l,p,pf',0,0)$,
$(2pk\pm j,2p,2pf',0,0)$. Applying Lemma \ref{lemma 0,1,c}
and \ref{lemma 1,2,c} we obtain that these vectors are isometric
respectively to $(0,1,r,0,0)$, $(1,2,2s,0,0)$, $(l,p,pt,0,0)$
$(j,2p,2pu,0,0)$.\erem
\begin{rem}{\rm The vector $(ts,t,f,0,0)$ is isometric to
$(0,t,*,0,0)$ by applying $s$-times $R_v\circ D$.}\end{rem}
\begin{prop}\label{prop: possible orbit}
Let $p$ be a prime number. The orbits of the following vectors of
$T_{2p}$
under isometries of $T_{2p}$ are all disjoint:\\
$\bullet$ $v_0=(1,0,0,0,0)$;\\
$\bullet$ $v_{1}=(0,1,r,0,0);$ \\ $\bullet$ $v_2=(1,2,2s,0,0)$;\\
$\bullet$ $v_p=(l,p,pt,0,0)$, where $0<l\leq
\llcorner p/2\lrcorner$;\\
$\bullet$ $v_{2p}=(j,2p,2pu,0,0)$, where $0<j< p$, $j\equiv 1\mod
2$.\end{prop} \bprf If two vectors $x$, $y$ of $T_{2p}$ are isometric, then
$x^2=y^2$ and the discriminants of the lattices orthogonal to $x$ and $y$ are equal: $d(x^{\perp})=d(y^{\perp})$. We resume the properties of the vectors $v_i$ 
in the following table: \small
\begin{eqnarray*}
\begin{array}{|l|l|l|l|l|l|}
\hline
v\!\!\!&v_0\!\!&v_1\!\!&v_2\!\!&v_p\!\!&v_{2p}\\
\hline
v^2\!\!\!\!&-2p\!\!&2r\!\!&-2p+8s\!\!&-2pl^2+2p^2t\!\!&-2pj^2+8p^2u\\
\hline
v^{\perp}\!\!\!\!&U\oplus U\!\!&\langle-2p\rangle\oplus \langle
-2r\rangle\oplus U\!\!&\left[\begin{array}{cc}
-2p\!\!&p\\p\!\!\!\!&-2s\end{array}\right]\oplus U
\!\!&\left[\begin{array}{cc} -2p\!\!&2l\\2l\!\!&-2t\end{array}\right]\oplus U \!\!&\left[\begin{array}{cc} -2p\!\!&j\\j\!\!&-2u\end{array}\right]\oplus U \\
\hline
d(v^{\perp})\!\!\!&1\!\!&-4pr\!\!&-p(4s-p)\!\!&-4(pt-l^2)\!\!&-4pu+j^2\\
\hline
\end{array}
\end{eqnarray*}
\normalsize For each copy of vectors $x$ and $y$ chosen from
$v_0$, $v_1$, $v_2$, $v_p$, $v_{2p}$ the conditions $x^2=y^2$ and
$d(x^{\perp})=d(y^{\perp})$ are incompatible. For example let us
analyze the case of $v_p$ and $v_{2p}$, the other cases are similar. We have:
$$
-2pl^2+2p^2t=-2pj^2+8p^2u\mbox{ and } \ 4(pt-l^2)=-4pu+j^2.
$$
By the first equation $-l^2+pt=-j^2+4pu$. Substituting in the
second equation we obtain $5(pt-l^2)=0$ and so $pt=l^2$. This
implies $p|l^2$ and so $p|l$. Since $l\leq \llcorner p/2\lrcorner$ this is
impossible.\erem

The previous results imply the following proposition:
\begin{prop}\label{prop: possible orbit of -4+U+U}
A primitive vector $(a_0,a_1,a_2,a_3,a_4)$ of the lattice $\langle -2p\rangle\oplus U\oplus U$ is isometric to
exactly one of the vectors:\\
$\bullet$ $v_0=(1,0,0,0,0)$;\\
$\bullet$ $v_{1}=(0,1,r,0,0)$ where $2a_1a_2+2a_3a_4-2pa_0^2=2r$;\\
$\bullet$ $v_2=(1,2,2s,0,0)$ where $2a_1a_2+2a_3a_4-2pa_0^2=-2p+8s$;\\
$\bullet$ $v_p=(l,p,pt,0,0)$, where $0<l\leq
\llcorner p/2\lrcorner$ and $2a_1a_2+2a_3a_4-2pa_0^2=-2pl^2+2p^2t$;\\
$\bullet$ $v_{2p}=(j,2p,2pu,0,0)$, where $0<j\leq p$ and
$2a_1a_2+2a_3a_4-2pa_0^2=-2pj^2+8p^2u$.
\end{prop}

\begin{rem}{\rm In particular in the case $p=2$ the only possibilities are
the vectors $(1,0,0,0,0)$, $(0,1,r,0,0)$, $(1,2,2s,0,0)$ and
$(1,4,4u,0,0)$.}\end{rem}
\subsection{The family} Let us denote by $\mathcal{L}_{r,w_i}^{2d}$ the overlattice of index $r$ of $\Z w_i\oplus\Omega_{(\Z/2\Z)^4}$, with $w_i^2=2d$ described in Theorem \ref{reticoli}. If $X$ is a K3 surface such that $NS(X)\simeq \mathcal{L}_{r,w_i}^{2d}$ for a certain $r=2,4,8$ and $i=1,2,3$, then $\Omega_{(\Z/2\Z)^4}$ is clearly primitively embedded in $NS(X)$ and thus $X$ admits $(\Z/2\Z)^4$ as group of symplectic automorphisms (cf. \cite[Theorem 4.15]{Nikulinsymplectic}). Hence, the lattices $\mathcal{L}_{r,w_i}^{2d}$ determine the family of algebraic K3 surfaces admitting a symplectic action of $(\Z/2\Z)^4$. More precisely: 
\begin{cor} \label{corollary: family X} The families of algebraic K3 surfaces admitting a symplectic action of $(\Z/2\Z)^4$ are the families of $\left(\mathcal{L}_{r,w_i}^{2d}\right)$-polarized K3 surfaces, for a certain $r=2,4,8$, $i=1,2,3$, $d\in 2\N_{>0}$. In particular the moduli space has a countable numbers of connected components of dimension 4.
\end{cor}
\begin{rem}\label{rem: quartic with Z/2Z4}{\rm If one fixes the value of $d$, then there is a finite number of possibilities for $r$ and $w_i$: for example if $d=2$, then $r=2$ and $i=1$, $w_1=(0,1,1,0,0,0,0)$. This implies that the family of quartic surfaces in $\mathbb{P}^3$ admitting a symplectic action of $(\Z/2\Z)^4$ has only one connected component of dimension 4. In \cite{E} the family of quartics invariant for the Heisenberg group($\simeq(\Z/2\Z)^4$) is described and since it is a 4-dimensional family of K3 surfaces admitting $(\Z/2\Z)^4$ as group of symplectic automorphisms we conclude that the family presented in \cite{E} is the family of the $(\mathcal{L}_{2,w_1}^4)$-polarized K3 surfaces.  The N\'eron--Severi group of such a K3 surfaces are generated by conics as proved in \cite[Corollary 7.4]{E}. 
}\end{rem}
\subsection{The subfamily of Kummer surfaces}
By Corollary \ref{cor: family Kummer}, for every non negative integer $d$ there exists a connected component of the moduli space of Kummer surfaces, which we called $\mathcal{F}_d$ and is the family of the $\mathcal{K}'_{4d}$-polarized K3 surfaces. For every $d$ the component $\mathcal{F}_d$ is 3-dimensional and by Proposition \ref{prop: every Kummer Z2Z4} it is contained in a connected component of the moduli space of K3 surfaces $X$ admitting $G$ as group of symplectic automorphisms. The following proposition identifies the components of the moduli space of K3 surfaces with a symplectic action of $G$ which contain $\mathcal{F}_d$:
\begin{prop}\label{prop: Kummer subfamily X}
The family of the $\mathcal{K}_{4d}'$-polarized Kummer surfaces is a codimension one subfamily of the following families: the $\left(\mathcal{L}_{2,w_1}^{4d}\right)$-polarized K3 surfaces; the $\left(\mathcal{L}_{4,w_2}^{8(2d-1)}\right)$-polarized K3 surfaces; the $\left(\mathcal{L}_{8,w_3}^{8(8d-1)}\right)$-polarized K3 surfaces. 
\end{prop}
\proof It suffices to show that there exists a primitive embedding $\mathcal{L}_{i,w_j}^{h}\subset\mathcal{K}_{4d}'$  or equivalently a primitive embedding $\left(\mathcal{K}_{4d}'\right)^{\perp}\subset\left(\mathcal{L}_{i,w_j}^{h}\right)^{\perp}$ for $(i,j,h)=(2,1,4d), (4,2,8(2d-1)), (8,3,8(8d-1))$. We recall that $\left(\mathcal{K}_{4d}'\right)^{\perp}\simeq \langle -4d\rangle\oplus U(2)\oplus U(2)$ and $\left(\mathcal{L}_{i,w_j}^{h}\right)^{\perp}$ is the transcendental lattice of the generic K3 surface $X$ described in Theorem \ref{reticoli}. With the notation of Theorem \ref{reticoli} sending a basis of 
$\langle -4d\rangle\oplus U(2)\oplus U(2)$ to the basis vectors $(0,1,0,0,0,0)$, $(0,0,1,0,0,0)$, $(0,0,0,1,0,0)$, $(0,0,0,0,1,0)$, $(0,0,0,0,0,1)$ of $\left(\mathcal{L}_{i,w_j}^{h}\right)^{\perp}$ with $t=d$, $s=d$, $u=d$ if $i=2,4,8$ respectively, we obtain an explicit primitive embedding of $\left(\mathcal{K}_{4d}'\right)^{\perp}$ in $\left(\mathcal{L}_{i,w_j}^{h}\right)^{\perp}$.
\eprf
We observe that the sublattice of $NS(Km(A))$ invariant for the action induced by the translation by the two torsion points on $A$, i.e., invariant for the action of $G$ defined in Proposition \ref{prop: every Kummer Z2Z4}, is generated by $H$ and $\frac{1}{2}(\sum_{p\in(\Z/2\Z)^4}K_p)$. Indeed $H$ is the image of the generator of $NS(A)$ by the map ${\pi_A}_*$, with the notation of diagram \eqref{diagkummer}. Thus, the lattice $\Omega_{(\Z/2\Z)^4}$ is isometric to $\langle H,\frac{1}{2}(\sum_{p\in(\Z/2\Z)^4K_p})\rangle^{\perp}\cap NS(Km(A))$ and in fact the lattice $\mathcal{L}_{2,w_1}^{2d}$ (which contains $\Omega_{(\Z/2\Z)^4}$ and an ample class) is isometric to $\langle\frac{1}{2}(\sum_{p\in(\Z/2\Z)^4}K_p)\rangle^{\perp}\cap NS(Km(A))$. 
\begin{rem}{\rm The previous proposition implies that the family of the Kummer surfaces of a $(1,d)$-polarized Abelian surface is contained in at least three distinct connected components of the family of K3 surfaces admitting a symplectic action of $G$. In particular, the intersection among the connected components of such family of K3 surfaces is non empty and of dimension 3.} 
\end{rem}
\subsection{Enriques involution}
In Section \ref{sec: automorphisms on Kummer} we have seen the result of Keum, \cite{Keum Enriques}: Every Kummer surface admits an Enriques involution. We now prove that this property holds more in general for the K3 surfaces admitting $(\Z/2\Z)^4$ as group of symplectic automorphisms and minimal Picard number.
\begin{theorem} Let $X$ be a K3 surface admitting $(\Z/2\Z)^4$ as group of symplectic automorphisms and such that $\rho(X)=16$, then $X$ admits an Enriques involution.\end{theorem}
\proof By Proposition \ref{prop: kummer enriques embedding} it suffices to prove that the transcendental lattice of $X$ admits a primitive embedding in $U\oplus U(2)\oplus E_8(-2)$ whose orthogonal does not contain vectors of length $-2$. The existence of this embedding can be proved as in \cite{Keum Enriques}. We briefly sketch the proof. Let $Q$ be one of the following lattices: $\langle -4\rangle\oplus \langle-2t\rangle$, $\left[\begin{array}{rr}-4&2\\2&-2s\end{array}\right]$, $\left[\begin{array}{rr}-4&1\\1&-2u\end{array}\right]$. The transcendental lattice of $X$ is $(U^2\oplus Q)(2)$. It suffices to prove that there exists a primitive embedding of $U(2)\oplus Q(2)$ in $U\oplus E_8(-2)$. The lattice $\langle -2\rangle\oplus Q$ is an even lattice with signature $(0,3)$. By \cite[Theorem 14.4]{Nikulinbilinear}, there exists a primitive embedding of $\langle -2\rangle\oplus Q$ in $E_8(-1)$, which induces a primitive embedding of $\langle -4\rangle\oplus Q(2)$ in $E_8(-2)$. Let $b_1$, $b_2$, $b_3$ be the basis of $\langle -4\rangle\oplus Q(2)$ in $E_8(-2)$. Let $e$ and $f$ be a standard basis of $U$ (i.e.\ $e^2=f^2=0$, $ef=1$). Then the vectors $e$, $e+2f+b_1$, $b_2$, $b_3$ give a primitive embedding of $U(2)\oplus Q(2)$ in $U\oplus E_8(-2)$ whose orthogonal complement does not contain vectors of length $-2$ (cf.\ \cite[\S 2, Proof of Theorem 2]{Keum Enriques}).\eprf

\section{The quotient K3 surface}\label{Y}
The surface $Y$ obtained as desingularization of the quotient
$X/(\Z/2\Z)^4$ contains 15 rational curves $M_i$, which are the
resolution of the 15 singular points of type $A_1$ on
$X/(\Z/2\Z)^4$. The minimal primitive sublattice of $NS(Y)$
containing these curves is denoted by $M_{(\Z/2\Z)^4}$. It is described in
\cite[Section 7]{Nikulinsymplectic} as an overlattice of the lattice
$\langle M_i\rangle_{i=1,\ldots, 15}$ of index $2^4$.

\begin{prop}\label{prop: overlattices of L+MG}
Let $Y$ be a K3 surface such that there exists a projective K3 surface $X$
and a symplectic action of $(\Z/2\Z)^4$ on $X$ with
$Y=\widetilde{X/(\Z/2\Z)^4}$. Then $\rho(Y)\geq 16$.\\
Moreover if $\rho(Y)=16$, let $L=M_{(\Z/2\Z)^4}^{\perp_{NS(Y)}}$. Then
$NS(Y)$ is an overlattice of index 2 of $\Z L\oplus M_{(\Z/2\Z)^4}$, where
$L^2=2d>0$. In particular, $NS(Y)$ is generated by $\Z L\oplus M_{(\Z/2\Z)^4}$
and by a class $(L/2,v/2)$, $v/2\in M_{(\Z/2\Z)^4}^{\vee}/M_{(\Z/2\Z)^4}$ (that is not trivial in $M_{(\Z/2\Z)^4}^{\vee}/M_{(\Z/2\Z)^4}$), $L^2\equiv -v^2\mod8$.
\end{prop}
\bprf A K3 surface $Y$ obtained as
desingularization of the quotient of a K3 surface $X$ by the symplectic
group of automorphisms $(\Z/2\Z)^4$, has $M_{(\Z/2\Z)^4}\subset NS(Y)$.
Since $M_{(\Z/2\Z)^4}$ is negative definite and $Y$ is projective (it is the quotient of $X$, which is
projective), there
is at least one class in $NS(Y)$ which is not in $M_{(\Z/2\Z)^4}$ so 
$\rho(Y)\geq
1+\rk M_{(\Z/2\Z)^4}=16$. In particular if $\rho(Y)=16$, then the
orthogonal complement of $M_{(\Z/2\Z)^4}$ in $NS(Y)$ is generated by a class with a positive
self intersection, hence $NS(Y)$ is either $\Z L\oplus M_{(\Z/2\Z)^4}$ or an
overlattice of $\Z L\oplus M_{(\Z/2\Z)^4}$ with a finite index. The
discriminant group of $M_{(\Z/2\Z)^4}$ is $(\Z/2\Z)^7$ by \cite[Section 7]{Nikulinsymplectic} and so the
discriminant group of the lattice $\Z L\oplus M_{(\Z/2\Z)^4}$ is
$(\Z/2d\Z)\oplus (\Z/2\Z)^7$. It has eight generators. If the 
lattice $\Z L\oplus M_{(\Z/2\Z)^4}$ is the N\'eron--Severi group of a K3 surface $Y$, then also the discriminant
group of $T_Y$ has eight generators, but $T_Y$ has rank $22-\rho(Y)=6$, so 
this is impossible. Hence $NS(Y)$ is an
overlattice of $\Z L\oplus M_{(\Z/2\Z)^4}$. The index of the inclusion and the costruction of the overlattice 
can be computed as in \cite[Proposition 2.1]{projectivemodels} or as in Theorem \ref{prop: mathcal K}.\eprf

The Kummer surfaces are also examples of K3 surfaces obtained as desingularization of the
quotient of K3 surfaces by the action of $(\Z/2\Z)^4$ as group 
of symplectic automorphisms, see Proposition \ref{prop: every Kummer Z2Z4}. 

In \cite[Sections 4.2, 4.3]{alikummer} the action of $G$ on the Kummer lattice and the construction of the surface $\widetilde{Km(A)/G}$ are described. The images of the curves $K_{a,b,c,d}$, $(a,b,c,d)\in(\Z/2\Z)^4$ on $Km(A)$ under the quotient map
$Km(A)\lra Km(A)/G$ is a single curve. This curve can be naturally identified with the curve $K_{0,0,0,0}$ on the minimal resolution $\widetilde{Km(A)/G}\cong Km(A)$ (see \cite{alikummer}). The minimal resolution contains also fifteen $(-2)$-curves coming from the blowing up of the nodes on $Km(A)/G$, which can be identified with $K_{e,f,g,h}$, $(e,f,g,h)\in(\Z/2\Z)^4\backslash\{(0,0,0,0)\}$. These are the fifteen $(-2)$-curves in $M_{(\Z/2\Z)^4}$, hence $M_{(\Z/2\Z)^4}=K_{(0,0,0,0)}^{\perp}\cap K$.
This identification allows us to identify the curves of
$M_{(\Z/2\Z)^4}$ with the points of the space
$(\Z/2\Z)^4\backslash\{(0,0,0,0)\}$, hence we denote them by $M_{a,b,c,d}$, $(a,b,c,d)\in
(\Z/2\Z)^4\backslash\{(0,0,0,0)\}$. More explicitly, we are identifying the curve $K_{a,b,c,d}$ with the curve $M_{a,b,c,d}$ for any $(a,b,c,d)\in(\Z/2\Z)^4\backslash \{(0,0,0,0)\}$. By \cite{Nikulinsymplectic}
the lattice $M_{(\Z/2\Z)^4}$ contains the 15 curves $M_{a,b,c,d}$,
$(a,b,c,d)\in (\Z/2\Z)^4\backslash\{(0,0,0,0)\}$, it is generated
by 11 of these curves and by 4 other classes which are linear combination of these curves with rational coefficients. These 4
classes have to be contained also in $K$ (because $M_{(\Z/2\Z)^4}\subset K$)
and hence they correspond to hyperplanes in $(\Z/2\Z)^4$ which do
not contains the point $(0,0,0,0)$ (because $K_{(0,0,0,0)}\not\in
M_{(\Z/2\Z)^4}$).\\
From now on $\bar{K}_W$ (resp. $\bar{M}_W$) denotes 
$\frac{1}{2}\sum_{p\in W}K_p$ (resp. $\frac{1}{2}\sum_{p\in W}M_p$) for a subset $W$ of $(\Z/2\Z)^4$
(resp. $W$ a subset of $(\Z/2\Z)^4\backslash\{(0,0,0,0)\}$). We determine the 
orbits of elements in the discriminant group of $M_{(\Z/2\Z)^4}$ and its isometries
using the ones of $K$. 

\begin{prop}\label{prop: orbits MG} With respect to the group of isometries of $M_{(\Z/2\Z)^4}$ there are exactly six distinct orbits in the discriminant group $M_{(\Z/2\Z)^4}^\vee/M_{(\Z/2\Z)^4}$.
\end{prop}
\proof 
Let $W$ be one of the following subspaces:
\begin{itemize}
\item[1)] $W=(\Z/2\Z)^4$; \item[2)] $W$ is a hyperplane in
$(\Z/2\Z)^4$; \item[3)] $W$ is a 2-dimensional plane in
$(\Z/2\Z)^4$; \item[4)] $W=V\ast V'$ where $V$ and $V'$ are
2-dimensional planes and $V\cap V'$ is a point.
\end{itemize}
By Remark \ref{rem: properties of K} the classes $\bar{K}_W$ are in $K^{\vee}$ and if $W$ is as in $1)$ or $2)$ the classes $\bar{K}_W\in K$, and thus they are trivial in $K^{\vee}/K$. If $W$ is such that $(0,0,0,0)\notin W$, then the class
$\bar{M}_W=\bar{K}_W$ is contained in $M_{(\Z/2\Z)^4}^{\vee}$. Indeed it is a linear
combination with rational coefficients of the curves $M_{(a,b,c,d)}$
with $(a,b,c,d)\in(\Z/2\Z)^4\backslash\{(0,0,0,0)\}$, i.e. it is in 
$M_{(\Z/2\Z)^4}\otimes \Q$. Moreover it has an integer intersection with
all the classes in $K$ and so in particular with all the classes
in $M_{(\Z/2\Z)^4}\subset K$, i.e. it is in $M_{(\Z/2\Z)^4}^{\vee}$. We observe that if
$W$ is a hyperplane (as in case $2)$) and it is such that
$(0,0,0,0)\notin W$, then the class $\bar{M}_W$ is a class in
$M_{(\Z/2\Z)^4}$ (and hence trivial in the discriminant group, see Remark \ref{rem: properties of K}).\\
If $(0,0,0,0)\in W$, let $W'$ be
$W'=W-\{(0,0,0,0)\}$. The class $\bar{M}_{W'}$ is a class in
$M_{(\Z/2\Z)^4}^{\vee}$. Indeed it is clear that $\bar{M}_{W'}\in M_{(\Z/2\Z)^4}\otimes
\Q$ has an integer intersection with all the classes
$M_{(a,b,c,d)}\in M_{(\Z/2\Z)^4}$, $(a,b,c,d)\in (\Z/2\Z)^4\backslash\{(0,0,0,0)\}$. 
Let $Z$ be a hyperplane of $(\Z/2\Z)^4$ which does not contain
$(0,0,0,0)$. Since $\bar{M}_Z\in M_{(\Z/2\Z)^4}$ we have to check that $\bar{M}_{W'}\cdot\bar{M}_Z\in\Z$.  We recall that
$\bar{K}_W$ is in $K^{\vee}$ and so it has an integer intersection
with all the classes $\bar{K}_Z$. This means that $W\cap Z$ is
made up of an even number of points. Since $(0,0,0,0)\notin Z$,
$(0,0,0,0)\notin W\cap Z$ and hence $W'\cap Z$ is an even number
of points. This implies that $\bar{M}_{W'}\cdot \bar{M}_Z\in \Z$.\\
If $\bar{M}_W\in M_{(\Z/2\Z)^4}^{\vee}$, hence either
$\bar{K}_W$ or $\bar{K}_{W\cup \{(0,0,0,0)\}}$ is in $K^{\vee}$.
Indeed by Remark \ref{rem: properties of K} the Kummer
lattice is generated by the curves $K_{(a,b,c,d)}$, $(a,b,c,d)\in(\Z/2\Z)^4$, by 4 classes
of type $\bar{K}_{W_i}$ where $W_i$ is the hyperplane $a_i=0$, $i=0,1,2,3$ (see the notation of Remark \ref{rem: properties of K}) and
by the class $\bar{K}_{(\Z/2\Z)^4}$. This is clearly equivalent
to say that $K$ is generated by the curves $K_{(a,b,c,d)}$, by 4
classes of type $\bar{K}_{W'_i}$ where $W'_i$ is the hyperplane
$a_i=1$ and by the class $\bar{K}_{(\Z/2\Z)^4}$. If $\bar{M}_W\in
M_{(\Z/2\Z)^4}^{\vee}$, then $\bar{M}_W\cdot \bar{M}_{W'_i}=\bar{K}_W\cdot
\bar{K}_{W'_i}\in \Z$.
Moreover, since $(0,0,0,0)\notin W'_i$, we have also
$\bar{K}_{W\cup \{(0,0,0,0)\}}\cdot \bar{K}_{W'_i}\in \Z$. To
conclude that either $\bar{K}_W$ or $\bar{K}_{W\cup
\{(0,0,0,0)\}}$ is in $K^{\vee}$, it suffices to prove either that
$\bar{K}_W\cdot \bar{K}_{(\Z/2\Z)^4}\in \Z$ or $\bar{K}_{W\cup
\{(0,0,0,0)\}}\cdot \bar{K}_{(\Z/2\Z)^4}\in \Z$. This is clear,
indeed $\bar{K}_W\cdot \bar{K}_{(\Z/2\Z)^4}\in \Z$ if and only if
$W$ consists of an even number of points. If it is not, clearly
$W\cup\{(0,0,0,0)\}$ consists of an even number of
points. Thus, the classes $\bar{M}_W$ are in $M_{(\Z/2\Z)^4}^{\vee}$ for the following
subspaces:

\begin{itemize} 
\item[1)] $W=(\Z/2\Z)^4\backslash\{(0,0,0,0)\}$; 
\item[2a)] $W$ is an
hyperplane in $(\Z/2\Z)^4$, $(0,0,0,0)\notin W$;
\item[2b)]
$W\backslash\{(0,0,0,0)\}$ where $W$ is a hyperplane in $(\Z/2\Z)^4$,
$(0,0,0,0)\in W$; 
\item[3a)] $W$ is a 2-dimensional plane in
$(\Z/2\Z)^4$ and $(0,0,0,0)\notin W$; 
\item[3b)] $W\backslash\{(0,0,0,0)\}$
where $W$ is a 2-dimensional plane in $(\Z/2\Z)^4$ and
$(0,0,0,0)\in W$; 
\item[4a)] $W=V\ast V'$ where $V$ and $V'$ are
2-dimensional planes and $V\cap V'$ is a point, $(0,0,0,0)\notin
V\ast V'$;
\item[4b)] $W\backslash\{(0,0,0,0)\}$ where $W=V\ast V'$, $V$ and
$V'$ are 2-dimensional planes and $V\cap V'$ is a point,
$(0,0,0,0)\in V\ast V'$.\end{itemize}
Each of these cases corresponds 
to a class of equivalence in the quotient $M_{(\Z/2\Z)^4}^{\vee}/M_{(\Z/2\Z)^4}$, here
we consider these equivalence classes. We will denote with $H$ an
hyperplane of $(\Z/2\Z)^4$ such that $(0,0,0,0)\notin H$. We observe that $\bar{M}_{W\ast H}\equiv \bar{M}_W+\bar{M}_H\mod \oplus_p\Z M_p$. Clearly the two classes $\bar{M}_W$ and $\bar{M}_{W\ast H}$ coincide in $M_{(\Z/2\Z)^4}^{\vee}/M_{(\Z/2\Z)^4}$ if $\bar{M}_H\in M_{(\Z/2\Z)^4}$. 
Let $n$ be the cardinality of $W\cap H$, $m$ be the number of curves $M_{(a,b,c,d)}$ appearing in $\bar{M}_{W\ast H}$ with a rational, non integer coefficient. In the following table we resume the classes of $M_{(\Z/2\Z)^4}^{\vee}$ which coincide modulo $M_{(\Z/2\Z)^4}$ and for each of them we give the value $discr$ of the discriminant form on it.
The first value of $m$ in the table is the number of curves in $\bar{M}_W$ and we put a $0$ for $n$; 
$$
\begin{array}{|c|c|c|c|}
\hline
\mbox{Case}&n&m&discr\\
\hline
1); 2b)&0;0,8;3&15; 7,7;7&\frac{1}{2}\\
\hline
3a)&0,4,2,0&4,4,8,12&0\\
\hline
3b)&0,0,2&3,7,11&\frac{1}{2}\\
\hline
4a)&4,2&6,10&1\\
\hline
4b)&4,2&5,9&-\frac{1}{2}\\
\hline
\end{array}$$ 
Indeed by Remark \ref{rem: properties of K} the orbit of elements in the discriminant
group $A_K$ of $K$ 
are three up to isometries. To prove the latter one considers the action of the group
$GL(4,\Z/2\Z)$ on $(\Z/2\Z)^4$, which in fact we can identify with a
subgroup of $O(A_K)$. Since
$GL(4,\Z/2\Z)$ fixes $(0,0,0,0)$, it acts also on
$(\Z/2\Z)^4\backslash\{(0,0,0,0)\}$ and so we can identify it with a subgroup of
$O(A_{M_{(\Z/2\Z)^4}})$. This means that under the action of $GL(4,\Z/2\Z)$ we have at
most six orbits, associated to the cases 1;2b), 2a), 3a), 3b),
4a), 4b). We observe that the orbit of 2a) is the one of class
$0\in A_{M_{(\Z/2\Z)^4}}$. We show now that all these orbits
are disjoint, so we have exactly 6 (5 non trivial) orbits in
$M_{(\Z/2\Z)^4}^{\vee}/M_{(\Z/2\Z)^4}$.
 One can check by a direct computation that the classes of the cases 1) and 2b) coincide in the quotient. The classes in $M_{(\Z/2\Z)^4}$ with self intersection $-2$ 
are only $\pm M_{(a,b,c,d)}$,
$(a,b,c,d)\in(\Z/2\Z)^4\backslash\{(0,0,0,0)\}$. Indeed each class in $M_{(\Z/2\Z)^4}$
is a linear combination
$D=\sum_{(a,b,c,d)\in(\Z/2\Z)^4-\{(0,0,0,0)\}}
\alpha_{(a,b,c,d)}M_{a,b,c,d}$ with $\alpha_{(a,b,c,d)}\in
\frac{1}{2}\Z$. The condition
$-2=D^2=-2\sum_{(a,b,c,d)}\alpha_{(a,b,c,d)}^2$ implies that either
there is one $\alpha_{(a,b,c,d)}=\pm 1$ and the others are zero,
or there are four $\alpha_{(a,b,c,d)}$ equal to $\pm \frac{1}{2}$
and the others are zero. Since there are no classes in $M_{(\Z/2\Z)^4}$
which are linear combination with rational coefficients of only four
classes, we have $D=\pm M_{(a,b,c,d)}$ for a certain
$(a,b,c,d)\in(\Z/2\Z)^4\backslash\{(0,0,0,0)\}$. Since the isometries of $M_{(\Z/2\Z)^4}$ preserve the intersection product, they send the classes of the curves
$M_{(a,b,c,d)}$ either to the class of a  curve or to the opposite
of the class of a curve. In particular, there are no isometries
of $M_{(\Z/2\Z)^4}$ which identify classes associated to the six
cases 1);2b), 2a), 3a), 3b), 4a), 4b), indeed in each class there is some linear
combination with non integer coefficients of a different number of curves $M_{(a,b,c,d)}$ .\eprf

\begin{theorem}\label{preciseNS}
Let $Y$ be a projective K3 surface such that there exists a K3 surface $X$
and a symplectic action of $(\Z/2\Z)^4$ on $X$ with
$Y=\widetilde{X/(\Z/2\Z)^4}$ and let $\rho(Y)= 16$.\\
Then $NS(Y)$ is
generated by $\Z L\oplus M_{(\Z/2\Z)^4}$ with $L^2=2d>0$, and by a class
$(L/2,v/2)$, $0\not=v/2\in M_{(\Z/2\Z)^4}^{\vee}/M_{(\Z/2\Z)^4}$ with $L^2\equiv -v^2\mod8$. Up to isometry there
are only the following possibilities:
\begin{itemize} \item[i)] if $d\equiv 1 \mod 4$, then $v/2=\bar{M}_{W}$ where  $$W=\{(1,1,0,1), (1,1,1,0), (1,1,1,1),(1,0,0,0), (0,1,0,0)\}$$  (case 4b) of proof of Proposition \ref{prop: orbits MG}); 
\item[ii)] if $d\equiv 2 \mod 4$, then $v/2=\bar{M}_{W}$
where $$W=\{(0,0,0,1), (0,0,1,0), (0,0,1,1), (1,0,0,0), (0,1,0,0),
(1,1,0,0)\}$$ 
(case 4a) of proof of Proposition \ref{prop: orbits MG});
\item[ ] if $d\equiv 3 \mod 4$, then: either\\
{\rm iii)} $v/2=\bar{M}_{W}$ where
$$W=\{(0,0,0,1), (0,0,1,0), (0,0,1,1)\}$$ 
(case 3b) proof of Proposition \ref{prop: orbits MG}), or\\ {\rm iv)} $v/2=\bar{M}_{W}$ where
$$W=(\Z/2\Z)^4-\{(0,0,0,0)\}$$ 
(case 1-2b)) proof of Proposition \ref{prop: orbits MG}); \item[v)] if $d\equiv 0 \mod 4$, then
$v/2=\bar{M}_{W}$ where 
$$W=\{(1,1,0,0), (1,1,1,0),
(1,1,0,1),(1,1,1,1)\}$$ (case 3a) of proof of Proposition \ref{prop: orbits MG}).
\end{itemize}
Moreover for each $d\in \N$ there exists a K3 surface $S$ such
that $NS(S)$ is an overlattice of index two of the lattice
$\langle 2d\rangle\oplus M_{(\Z/2\Z)^4}$.\\
In cases {\rm i), ii), iii), v)}, $T_Y\simeq U(2)\oplus U(2)\oplus \langle -2\rangle\oplus\langle -2d\rangle $.
In case {\rm iv)} denote by $q_2$ the discriminant form of $U(2)$ then the discriminant group of $T_Y$ is $(\Z/2\Z)^5\oplus \Z/2d\Z$ with discriminant form $q_2\oplus q_2\oplus \left(\begin{array}{cc} 0&1/2\\1/2&(-d-1)/2d\end{array}\right)$.
\end{theorem}

\bprf In Proposition \ref{prop: overlattices of L+MG} we proved
that the lattice $NS(Y)$ has to be an overlattice of index 2 of
$\Z L\oplus M_{(\Z/2\Z)^4}$. The unicity of the choice of $v$ up to isometry
depends on the description of the orbit of the group of the
isometries of $M_{(\Z/2\Z)^4}^{\vee}/M_{(\Z/2\Z)^4}$ given in Proposition \ref{prop: orbits MG}.\\
By an explicit computation one can show that the discriminant group of the overlattices described in 
 i), ii), iii), iv), v) is $(\Z/2\Z)^5\oplus (\Z/2d\Z)$ and the discriminant form
in all the cases except iv) is $q_2\oplus q_2\oplus 1/2 \oplus 1/2d$. In the case iv) the discriminant  
form is those described in the statement. In any case by \cite[Theorem 1.14.4 and Remark 1.14.5]{Nikulinbilinear}
the overlattices have a unique primitive embedding in the K3 lattice $\Lambda_{K3}$, hence by the surjectivity of the
period map there exists a K3 surface $S$ as in the statement of the theorem. Moreover by \cite[Theorem 1.13.2 and 1.14.2]{Nikulinbilinear} the transcendental lattice is uniquely determined by signature and discriminant form. This concludes the proof. \eprf 

\begin{rem}{\rm The Kummer surfaces appear as specializations of the surfaces $Y$ as in Proposition \ref{preciseNS} such that $d\equiv 0\mod 2$. Indeed, let us consider the surface $Y$ such that $d=2d'$. The transcendental lattice of a generic Kummer of a $(1,d')$-polarized Abelian surface is $T_{Km(A)}\simeq U(2)\oplus U(2)\oplus \langle -4d'\rangle$, and it is clearly primitively embedded in  $T_Y\simeq U(2)\oplus U(2)\oplus \langle -2\rangle \oplus \langle -4d'\rangle$.}\end{rem}

\subsection{Ampleness properties}
As in Section \ref{ampleness}, we can prove that certain divisors on $Y$ are ample (or nef or nef and big) using the description of the N\'eron--Severi group of $Y$ given in Theorem \ref{preciseNS}. The ample (or nef or nef and big) divisors define projective models, which can be described in the same way as in Section \ref{examplespic17}, where we described projective models of the Kummer surfaces.

\begin{prop}\label{prop: ampleness on Y} With the notation of Theorem \ref{preciseNS}, the following properties for divisors on $Y$ hold:
\begin{itemize}\item $L$ is pseudo ample and it has no fixed components; \item the divisor $D:=L-(M_1+\ldots + M_r)$, $1\leq r\leq 15$ is pseudo ample if $d>r$; \item let $\bar{D}:=\left(L-M_1-\ldots -M_r\right)/2\in NS(Y)\otimes\Q$; if $\bar{D}\in NS(Y)$, then it is pseudo ample if $d>r$.\end{itemize}
\end{prop}

\subsection{K3 surfaces with 15 nodes}
Here we show that a K3 surface with 15 nodes (resp. with 15 disjoint irreducible rational curves) is in fact the quotient (resp. the desigularization of the quotient) of a K3 surface by a symplectic action of $(\Z/2\Z)^4$. This is in a certain sense the generalization of a similar result for Kummer surface (cf. Section \ref{subsec:  K3 with 16 nodes}).

\begin{theorem}\label{15nodi}
Let $Y$ be a projective K3 surface with 15 disjoint smooth rational curves $M_i$, $i=1,\ldots, 15$ or equivalently a K3 surface admitting  a singular model with $15$ nodes. 
Then:\\
{\bf 1)} $NS(Y)$ contains the lattice $M_{(\Z/2\Z)^4}$.\\
{\bf 2)} There exists a K3 surface $X$ with a $G=(\mathbb{Z}/2\mathbb{Z})^4$ symplectic action, such that $Y$ is the minimal resolution of the quotient  $X/G$.\\
\end{theorem}
\bprf
{\bf 1)} Let $Q$ be the orthogonal complement in $NS(Y)$ to $\oplus_{i=1}^{15}\Z M_i$  and $R$ be the lattice $Q\oplus\left( \oplus_{i=1}^{15} \Z M_i\right)$. Observe that $NS(Y)$ is an overlattice of finite index of $R$ and $R^{\vee}/R\cong Q^{\vee}/Q\oplus(\Z/2\Z)^{\oplus 15}$ so $l(R)=l(Q)+15$. Let us denote by $k$ the index of $R$ in $NS(Y)$, thus $l(NS(Y))=l(Q)+15-2k$. On the other hand the rank of the transcendental lattice is $22-\rk(R)=7-\rk(Q)$. Hence $l(Q)+15-2k\leq 7-\rk(Q)$. Thus $k\geq \left(8+l(Q)+\rk(Q)\right)/2$.  We observe that $k$ is the minimum number of divisible class we have to add to $R$ in order to obtain $NS(Y)$. By definition the lattice $Q$ is primitive in $NS(Y)$, thus the non trivial classes that we can add to $R$ in order to obtain overlattices are either classes in $(\oplus_i \Z M_i)^{\vee}/(\oplus_i \Z M_i)$ or classes like $v+v'$, where $v'\in Q^{\vee}/Q$ and $v\in (\oplus_i \Z M_i)^{\vee}/(\oplus_i \Z M_i)$ is non trivial. By construction the independent classes of the second type are at most $l(Q)$ and thus there are at least $ (\left(8+l(Q)+\rk(Q)\right)/2)-l(Q)=\left(8+\rk(Q)-l(Q)\right)/2$ classes which are in $(\oplus_i \Z M_i)^{\vee}/(\oplus_i \Z M_i)$. We recall that $\rk(Q)-l(Q)\geq 0$ and hence there are at least 4 classes which are rational linear combinations of the curves $M_i$. By \cite[Lemma 3]{nikulinKummer} such a class in  $NS(Y)$ can only contain $16$ or $8$ classes. Since $16$ is not possible in this case, all these classes contain eight $(-2)$-curves.  Let $u_j$, $j=1,2,3,4$ be 4 independent classes in $\left(\oplus_i \Z M_i\right)^{\vee}/ \left(\oplus \Z M_i\right)$ such that the $u_j$ are contained in $NS(Y)$. For each $j\neq h$, $j,h=1,2,3,4$, there are exactly 4 rational curves which are summands of both $u_i$ and $u_j$, otherwise the sum $u_i+u_j\in NS(Y)$ contains half the sum of $k'$ disjoint rational curves for $k'\neq 8$, which is absurd. It is now a trivial computation to show that there are at most 4 independent classes (and thus exactly 4) as required in $\left(\oplus_i \Z M_i\right)^{\vee}/ \left(\oplus \Z M_i\right)$ and that for each choice of these 4 classes $u_i$, the lattice obtained adding the classes $u_i$, $i=1,2,3,4$  to $\oplus_i \Z M_i$ is exactly $M_{(\Z/2\Z)^4}$: indeed without loss of generality the first class can be chosen to be $u_1=\sum_{i=1}^8(M_i/2)$, thus the second class can be chosen to be $u_2=\sum_{i=1}^4 (M_i/2)+\sum_{j=9}^{12}(M_j/2)$. The third class has 4 curves in common with $u_1$ and with $u_2$ and thus can be chosen to be $u_3=(M_1+M_2+M_5+M_6+M_9+M_{10}+M_{13}+M_{14})/2$. Similarly, one determines the class $u_4=(M_1+M_3+M_5+M_7+M_9+M_{11}+M_{13}+M_{15})/2$.\\
{\bf 2)} We consider the double cover $\pi_1: Z_1\lra Y$ ramified on $2u_1$. Since $2u_1$ contains $8$ disjoint rational curves, $Z_1$ is smooth. Moreover the pullback $E_i$ of the curves  $M_i$, $i=1,\ldots,8$ have self intersection $-1$, hence these can be contracted to smooth points on a variety $Y_1$, and the covering involution that determines $\pi_1$ descends to a symplectic involution $\iota_1$ on $Y_1$ with $8$ isolated fixed points (cf. \cite[\S 3]{morrison}). The  divisors $2u_i$, $i=2,3,4$ contain each $4$ curves which are also in the support of $2u_1$.
We study the pull back of $2u_2$, for the other classes the study is similar. We have $2\pi_1^*(u_2)=\pi_1^*(2u_2) =2(E_1+\ldots+E_4)+M_5^1+M_5^2+M_6^1+M_6^2 +M_7^1+M_7^2 +M_8^1+M_8^2$, where $\pi_1(M^i_j)=M_j$ for $i=1,2$ and $j=5,6,7,8$. Hence the divisor  $M_5^1+M_5^2+M_6^1+M_6^2 +M_7^1+M_7^2 +M_8^1+M_8^2$ is divisible by $2$ in $NS(Z_1)$ and so its image is divisible by $2$ on $NS(Y_1)$. Doing the same construction as before, using this class we get a K3 surface $Y_2$ with an  action by a symplectic involution $\iota_2$. Observe that $\iota_1$ preserves the divisor  $M_5^1+M_5^2+M_6^1+M_6^2 +M_7^1+M_7^2 +M_8^1+M_8^2$ and so $\iota_1$ and $\iota_2$ commute on $NS(Y_2)$. Considering now the pull-back of $2u_3$ and $2u_4$ on $Y_2$ one can repeat the construction arriving at a K3 surface $X:=Y_4$ with an action by $(\Z/2\Z)^4$ and such that the quotient is $Y$. We observe that the smooth model of a K3 surface admitting a singular model with 15 nodes contains 15 disjoint rational curves and we proved that such a K3 surface is a $(\Z/2\Z)^4$ quotient of a K3 surface.  
\eprf

\begin{rem}{\rm Assume now that a K3 surface $S$ either has a lattice isometric to $M_{(\Z/2\Z)^4}$ primitively embedded in the N\'eron--Severi group or its N\'eron--Severi group is an overlattice of $Q\oplus \langle-2\rangle ^{15}$ for a certain lattice $Q$. Then the Theorems \ref{preciseNS},  \ref{15nodi} do not imply that $S$ is a $(\Z/2\Z)^4$ quotient of a K3 surface. Indeed in the proof of Theorem \ref{15nodi}, part 2), we used that the lattice $\langle -2\rangle^{15}$ (contained with index $2^4$ in $M_{(\Z/2\Z)^4}$) is generated by irreducible rational curves. In other words the description of the N\'eron--Severi group from a lattice theoretic point of view is not enough to obtain our geometric characterization. Thus we can not conclude that the family of the K3 surfaces which are $(\Z/2\Z)^4$ quotients of K3 surfaces coincides with the family of the K3 surfaces polarized with certain lattices.}\end{rem}

\begin{rem}{\rm In the proof of Theorem \ref{15nodi} we proved that if a K3 surface contains 15 disjoint rational curves $M_i$, then there are 15 subsets $\mathcal{S}_i$, $i=1,\ldots, 15$ of 8 of these curves which form an even set. Similarly if a K3 surface has 15 nodes there are 15 subsets of 8 of these nodes which form an even set. In \cite{barthevenset} and \cite{projectivemodels} some geometric properties of the even set of curves and nodes on K3 surfaces are described. For example if a quartic in $\mathbb{P}^3$ contains 8 nodes which form an even set, then the eight nodes are contained in an elliptic curve and there are three quadrics in $\mathbb{P}^3$ containing these nodes. Hence if a quartic in $\mathbb{P}^3$ has 15 nodes, each even set $\mathcal{S}_i$ has the previous properties.}\end{rem}
\begin{cor}\label{cor: 14 nodes}
Let $Y$ be a projective K3 surface with 14 disjoint smooth rational curves $M_i$, $i=1,\ldots,14$. Then:\\
{\bf 1)} $NS(Y)$ contains the lattice $M_{(\Z/2\Z)^3}$ which is the minimal primitive sublattice of the 
K3 lattice $\Lambda_{K3}$ that contains the 14 rational curves.\\
{\bf 2)} There exists a K3 surface with a $(\Z/2\Z)^3$ symplectic action, such that $Y$ is the minimal resolution of the quotient of $X$ by this group.      
\end{cor}
\bprf
The lattice  $M_{(\Z/2	\Z)^3}$  is described in \cite[Section 7]{Nikulinsymplectic}. The proof of {\bf 1)} and {\bf 2)} is essentially the same
as the proof of {\bf 1)} and {\bf 2)} of Theorem \ref{15nodi}.
\eprf
\begin{rem}{\rm  The result analogous to the one of Proposition \ref{15nodi} and Corollary \ref{cor: 14 nodes} does not hold considering 8 (resp. 12) disjoint rational curves, i.e. considering the group $\Z/2\Z$ (resp. $(\Z/2\Z)^2$):\\
{\bf 1)}  If a K3 surface is the minimal resolution of the quotient of a K3 surface by the group $\Z/2\Z$, then it admits a set of 8 disjoint rational curves but if a K3 surface admits a set of 8 disjoint rational curves, then it is not necessarily the quotient of a K3 surface by the group $\Z/2\Z$ acting symplectically. An example is given by the K3 surface with an elliptic fibration with 8 fibers of type $I_2$ and trivial Mordell--Weil group (cf. \cite[Table 1, Case 99]{shimada}): it contains 8 disjoint rational curves (a component for each reducible fibers), which are not an even set (otherwise the fibration admits a 2-torsion section).\\
{\bf 2)} If a K3 surface is the minimal resolution of the quotient of a K3 surface by the group $(\Z/2\Z)^2$, then it admits a set of 12 disjoint rational curves but if a K3 surface admits a set of 12 disjoint rational curves, then it is not necessarily the quotient of a K3 surface by a group $(\Z/2\Z)^2$ acting symplectically. Anyway it is surely the quotient of a K3 surface by $\Z/2\Z$ (the proof is again similar to the one of Theorem \ref{15nodi}). An example is given by the elliptic K3 surface with singular fibers $2 I_0^*+4 I_2$ which is the number 466 in Shimada's list, \cite{shimada}.  The components of multiplicity 1 of the fibers of type $I_0^*$ and a component for each fiber of type $I_2$ are 12 disjoint rational curves. There is exactly one set of 8 of these curves which is a 2-divisible class (the sum of the components of the $I_0^*$ fibers of multiplicity one). By using the Shioda-Tate formula one can easily show that there are no more divisible classes and hence the surface can not be the quotient of a K3 surface by $(\Z/2\Z)^2$. }\end{rem}


\section{The maps $\pi_*$ and $\pi^*$}\label{sec: maps}
In the previous two sections we described the family of the K3 surfaces $X$ admitting a symplectic action of $(\Z/2\Z)^4$ and the family of the K3 surfaces $Y$ obtained as desingularizations of the quotients of K3 surfaces by the group $(\Z/2\Z)^4$.  Here we explicitly describe the relation among these two families. More precisely in Section \ref{sec: K3 with Z2Z4 in general}, we described the quotient map $\pi:\widetilde{X}\ra Y$, which of course induces the maps $\pi_*$ and $\pi^*$ among the cohomology groups of the surfaces: here we describe these maps (similar results can be found in \cite{vGS} if the map $\pi$ is the quotient map by a symplectic involution).  With the notation of diagram  \eqref{diagrammone} we have:
\begin{prop}
The map $\pi_*:H^2(\widetilde{X},\Z)\lra H^2(Y,\Z)$  is induced by the map
$$
\begin{array}{ccc}
\pi_*:<-2>^{\oplus 16}\oplus U(2)^{\oplus 3}\oplus \left(<-1>^{\oplus 8}\right)^{\oplus 15}&\lra&<-2>\oplus U(32)^{\oplus 3}\oplus <-2>^{\oplus 15}\\
\pi_*:(k_1,\ldots, k_{16},u,\{n_{1j}\}_{1\leq j\leq 8},\ldots,\{n_{15j}\}_{1\leq j\leq 8})&\mapsto&(k,u,m_1,\ldots,m_{15})\\
\end{array}
$$
where $\pi_*(k_i)=k$, for all $i=1,\ldots, 16$; $\pi_*(n_{ij})=m_i$ for all $j=1,\ldots,8$,  $i=1,\ldots, 15$.\\
The map $\pi^*:H^2(Y,\Z)\lra H^2(\widetilde{X},\Z)$ is induced by the map
$$
\begin{array}{ccc}
\pi^*:<-2>\oplus U(32)^{\oplus 3}\oplus <-2>^{\oplus 15}&\hookrightarrow &<-2>^{\oplus 16}\oplus U(2)^{\oplus 3}\oplus \left(<-1>^{\oplus 8}\right)^{\oplus 15}\\
\pi^*:(k,u,m_1,\ldots,m_{15})&\mapsto& (k_1=k,\ldots,k_{16}=k, 16 u, \sum_{j=1}^{8}2n_{1j},\ldots,\sum_{j=1}^{8}2n_{15j})
\end{array}
$$
\end{prop}
\bprf
By \cite[Theorem 4.7]{Nikulinsymplectic} the action of $G$ on $\Lambda_{K3}$ does not depend on the K3 surface we have chosen, hence we can consider $X=Km(A)$  and $G$ realized as in Section \ref{sec: automorphisms on Kummer} (i.e. it is induced on $Km(A)$ by the translation by the 2-torsion points of the Abelian surface $A$).\\
{\bf 1.} $\pi_*$. We have seen that $G$ leaves $U(2)^{\oplus 3}$ invariant and in fact $H^2(X,\Z)^G\supset U(2)^{\oplus 3}$, however the map $\pi_*$ multiply the intersection form by 16. In fact for $x_1,x_2\in U(2)^{\oplus 3}$ we have:
$$
\pi^*\pi_*(x_1)=16x_1
$$
so using the projection formula
$$
(\pi_* x_1,\pi_*x_2)_Y=(\pi^*\pi_*x_1,x_2)_{\tilde{X}}=16(x_1,x_2).
$$
Since by taking  $X=Km(A)$ the classes in $<-2>^{\oplus 16}$ correspond to classes permuted by $G$ their image by $\pi_*$ is a single $(-2)$-class in $H^2(Y,\Z)$. Finally, the $120$ $(-1)$-classes which are the blow up of the points with a non trivial stabilizer on $X$ are divided in orbits of length eight and mapped to the same curve $m_i$ on $Y$. By using the projection formula and the fact that the stabilizer group of a curve $n_{ij}$ has order 2, we have
$$
(m_i,m_i)_Y=(\pi_*(n_{ij}), \pi_*(n_{ij}))_Y=(\pi^*\pi_*(n_{ij}),n_{ij})_{\tilde{X}}=(2(n_{i1}+\ldots +n_{i8}),n_{ij})_{\tilde{X}}=-2.
$$
{\bf 2.} $\pi^*$. Let $x\in U(32)^{\oplus 3}$ and $y\in U(2)^{\oplus 3}$ then
$$
(\pi^*x,y)_{\tilde{X}}=(x,\pi_*y)_Y=(x,y)_Y=16(x,y)_{\tilde{X}}
$$
so $\pi^*(x)=16 x$. Then we have $\pi^*(u)=16 u$ since $u$ is not
a class in the branch locus.
Finally
$$
(\pi^*(m_i),n_{hj})_{\tilde{X}}=(m_i,\pi_*(n_{hj}))_Y=(m_i,m_h)_Y=-2\delta_{ih}
$$
and $(\pi^*(m_i),k)_{\tilde{X}}=(\pi^*(m_i),u)_{\tilde{X}}=0$ for $u\in U(32)^{\oplus 3}$. 
Hence $\pi^*(m_i)$ is given as in the statement. \eprf 
\begin{rem}
{\rm The lattice
$R:=<-2>^{\oplus 16}\oplus U(32)^{\oplus 3}$ (which is an overlattice of
index $2^5$ of $K\oplus U(32)^{\oplus 3}$) has index $2^{23}$ in
$\Lambda_{K3}$. Here we want to consider the divisible classes
that we have to add to $<-2>^{\oplus 16}\oplus U(32)^{\oplus 3}$ to obtain
the lattice
$\Lambda_{K3}$. Consider the $\Z$ basis
$\{\omega_{ij}\}_{i\neq j}$ of $U(2)^3$ in $H^2(Km(A),\Z)$. Recall that we have
an exact sequance $0\ra A[2]\ra A\stackrel{\cdot 2}{\ra}A\ra 0$, which corresponds to the multiplication by 2 on each real coordinates of $A$. Thus, the copy of $U(32)^{\oplus 3}\subset H^2(Km(A/A[2]),\Z)$ is generated by
 $4\omega_{ij}$. Hence let $e_i,f_i$,
$i=1,2,3$ be the standard basis of each copy of $U(32)$, then the
elements:
$$
e_i/4,f_i/4
$$
are contained in $H^2(Y,\Z)$. Adding these classes to $R$ we find
$<-2>^{\oplus 16}\oplus U(2)^3$ as overlattice of index
$2^{12}$
of $R$.\\
In Remark \ref{rem: overlattice of -2^16+U(2)} the construction of 
the even unimodular overlattice $\Lambda_{K3}$ of $\langle
-2\rangle^{\oplus 16}\oplus U(2)^{\oplus 3}$ is described (we observe that
the index is $2^{11}$). In conclusion  we can construct explicitly the
overlattice $\Lambda_{K3}$ of $R$ and
extend the maps, $\pi_*$, $\pi^*$ to this lattice.}
\end{rem}

\section{Some explicit examples}\label{examples}
In this Section we provide geometrical examples of K3 surfaces $X$ with Picard number 16 admitting a symplectic action of $G=(\Z/2\Z)^4$ and of the quotient $X/G$, whose desingularization is $Y$.  We follow the notation of diagram \eqref{diagrammone} and  we denote by $L$ the polarization on $X$ orthogonal to the lattice $\Omega_{(\Z/2\Z)^4}$ and by $M$ the polarization on $Y$ orthogonal to the lattice $M_{(\Z/2\Z)^4}$.
\subsection{The polarization $L^2=4$, $M^2=L^2$.} We consider the projective space
$\PP^3$ and the group of transformations generated by:\\
\begin{eqnarray*}
\begin{array}{l}
(x_0:x_1:x_2:x_3)\mapsto (x_0:-x_1:x_2:-x_3)\\
(x_0:x_1:x_2:x_3)\mapsto (x_0:-x_1:-x_2:x_3)\\
(x_0:x_1:x_2:x_3)\mapsto (x_1:x_0:x_3:x_2)\\
(x_0:x_1:x_2:x_3)\mapsto (x_3:x_2:x_1:x_0)
\end{array}
\end{eqnarray*}
these transformations generate a group isomorphic to
$G=(\Z/2\Z)^4$. The invariant polynomials are
\begin{eqnarray*}\begin{array}{l}
p_0=x_0^4+x_1^4+x_2^4+x_3^4\\
p_1=x_0^2x_1^2+x_2^2x_3^2\\
p_2=x_0^2x_2^2+x_1^2x_3^2\\
p_3=x_0^2x_3^2+x_1^2x_2^2\\
p_4=x_0x_1x_2x_3
\end{array}
\end{eqnarray*}
Hence the generic $G$-invariant quartic K3 surface is a linear
combination:
$$
a_0(x_0^4+x_1^4+x_2^4+x_3^4)+a_1(x_0^2x_1^2+x_2^2x_3^2)+a_2(x_0^2x_2^2+x_1^2x_3^2)+a_3(x_0^2x_3^2+x_1^2x_2^2)+a_4x_0x_1x_2x_3=0.
$$
Since the only automorphism commuting with all the elements of the
group $G$ is the identity, the number of parameters in the equation
is $4$, which is also the dimension of the moduli space of the K3
surfaces with symplectic automorphism group $G$ and polarization
$L$ with $L^2=4$.\\ We study now the quotient surface. Observe
that the quotient of $\PP^3$ by $G$ is the Igusa quartic (cf.
\cite[Section 3.3]{hunt}), which is an order four relation between
the $p_i$'s, this is:
\begin{eqnarray*}
\begin{array}{ll}
\mathcal{I}_4\,: &16p_4^4+p_0^2p_4^2+p_1^2p_2^2+p_1^2p_3^2+p_2^2p_3^2-4(p_1^2+p_2^2+p_3^2)p_4^2-p_0p_1p_2p_3=0\\
\end{array}
\end{eqnarray*}
Hence the quotient is a quartic K3 surface which is a section of
the Igusa quartic by the hyperplane:
\begin{eqnarray*}
a_0p_0+a_1p_1+a_2p_2+a_3p_3+a_4p_4=0.
\end{eqnarray*}\\
The quartics in $\mathbb{P}^3$ admitting $(\Z/2\Z)^4$ as symplectic group of automorphisms are described in a very detailed way in \cite{E} (cf. also Remark \ref{rem: quartic with Z/2Z4}). We 
observe that the subfamily with $a_0=0$ is also a subfamily of the family of quartics considered by Keum
in \cite[Example 3.3]{Keum}. On this subfamily it is easy to identify an Enriques involution: this is the  
 standard Cremona transformation $(x_0:x_1:x_2:x_3)\lra (1/x_0:1/x_1:1/x_2:1/x_3)$. 
\subsection{The polarization $L^2=8$, $M^2=L^2/4=2$.}

Let $X$ be a K3 surface with a symplectic action of $G$ and $L^2=8$. There are two connected components of the moduli space of K3 surfaces with these properties (cf. Theorem \ref{reticoli} and Corollary \ref{corollary: family X}). One of them is realized as follows. 
Let us consider the complete intersection of three quadrics in $\PP^5$:
\begin{eqnarray*}
\left\{
\begin{array}{l}
\sum_{i=0}^5a_ix_i^2=0\\
\sum_{i=0}^5b_ix_i^2=0\\
\sum_{i=0}^5c_ix_i^2=0.
\end{array}
\right.
\end{eqnarray*}
with complex parameters $a_i$, $b_i$, $c_i$, $i=0,\ldots,5$. The group $G$ is realized as the transformations of $\PP^5$ changing an even number of signs in the coordinates.
To compute the dimension of the moduli space of these K3 surfaces we
must choose three independent quadrics in a six-dimensional 
space. Hence we must compute the dimension of the Grassmannian of
the subspaces of dimension three in a space of dimension six. This
is $3(6-3)=9$. Now the automorphisms of $\PP^5$ commuting with the
automorphisms generating $G$ are the diagonal $6\times
6$-matrices, hence we find the dimension $9-(6-1)=4$ as expected.\\
To determine the quotient, one sees that the invariant polynomials
under the action of $G$ are exactly the polynomials
$z_0^2,z_1^2,z_2^2,z_3^2,z_4^2,z_5^2$ and the product
$z_0z_1z_2z_3z_4z_5$. Denote them by $y_0,\ldots,y_5,t$ then there
is a relation
\begin{eqnarray*}
t^2=\prod_{i=0}^{5}y_i,
\end{eqnarray*}
and so we obtain a K3 surface which is the double cover of the
plane given by the intersection of the planes of $\PP^5$:
\begin{eqnarray*}
\left\{
\begin{array}{l}
\sum_{i=0}^5a_iy_i=0\\
\sum_{i=0}^5b_iy_i=0\\
\sum_{i=0}^5c_iy_i=0.
\end{array}
\right.
\end{eqnarray*}
The branch locus are six lines meeting at $15$ points, whose preimages under the double cover are the 15 nodes of the K3 surface. \\ 
We get a special subfamily of K3 surfaces considering as in Section \ref{jaco} a curve $\Gamma$ of genus 
2 with equation: 
\begin{eqnarray*}
y^2=\prod_{i=0}^{5}(x-s_i)
\end{eqnarray*}
with $s_i\in\mathbb{C}$, $s_i\not=s_j$ for $i\not= j$. This determines a family of Kummer surfaces with
$(\Z/2\Z)^4$ action and equations in $\PP^5$:
\begin{eqnarray*}
\left\{
\begin{array}{l}
z_0^2+z_1^2+z_2^2+z_3^2+z_4^2+z_5^2=0\\
s_0z_0^2+s_1z_1^2+s_2z_2^2+s_3z_3^2+s_4z_4^2+s_5z_5^2=0\\
s_0^2z_0^2+s_1^2z_1^2+s_2^2z_2^2+s_3^2z_3^2+s_4^2z_4^2+s_5^2z_5^2=0.
\end{array}
\right.
\end{eqnarray*}
The quotient surface also specializes to the double cover $t^2=\prod_i y_i$ of the
plane obtained as the intersection of the planes of $\PP^5$:
\begin{eqnarray*}
\left\{
\begin{array}{l}
y_0+y_1+y_2+y_3+y_4+y_5=0\\
s_0y_0+s_1y_1+s_2y_2+s_3y_3+s_4y_4+s_5y_5=0\\
s_0^2y_0+s_1^2y_1+s_2^2y_2+s_3^2y_3+s_4^2y_4+s_5^2y_5=0.\\
\end{array}
\right.
\end{eqnarray*}
As before the branch locus are 6 lines meeting at $15$ points, but in this case there is a conic tangent to the 6 lines. 

\addcontentsline{toc}{section}{ \hspace{0.5ex} References}

\end{document}